\documentclass[oneside, 12pt]{amsart}
\usepackage{amscd, amssymb, amsmath, mathrsfs}
\usepackage[english]{babel}
\usepackage{url}
\usepackage{tikz-cd}
\usepackage[pdftex, colorlinks=true,  citecolor=blue, linkcolor=blue, linktocpage=true ]{hyperref}

\newcommand\mylabel[1]{\label{#1}\marginpar{\raggedright\vspace{-1ex}\medskip\medskip\footnotesize \tt #1}}
\renewcommand\mylabel[1]{\label{#1}}

\newcommand{\mydate}{
\number\day\space
\ifcase\month \or January\or February\or March\or April\or May\or June\or July\or August\or September\or October\or November\or December\fi 
\space\number\year}

\DeclareUrlCommand\arXiv{\urlstyle{same}}

\newtheorem{theorem}{Theorem}[section]
\newtheorem*{maintheorem}{Theorem}
\newtheorem{lemma}[theorem]{Lemma}
\newtheorem{proposition}[theorem]{Proposition}
\newtheorem{corollary}[theorem]{Corollary}

\theoremstyle{definition}
\newtheorem{definition}[theorem]{Definition}

\newtheorem*{acknowledgement}{Acknowledgement}

\theoremstyle{remark}

\newcommand{\ZZ}{\mathbb{Z}}
\newcommand{\QQ}{\mathbb{Q}}

\newcommand{\CC}{\mathbb{C}}
\newcommand{\FF}{\mathbb{F}}

\newcommand{\PP}{\mathbb{P}}

\newcommand{\GG}{\mathbb{G}}

\newcommand{\ideala}{\mathfrak{a}}

\newcommand{\shA}{\mathscr{A}}

\newcommand{\shF}{\mathscr{F}}
\newcommand{\shG}{\mathscr{G}}

\newcommand{\shI}{\mathscr{I}}

\newcommand{\shM}{\mathscr{M}}

\newcommand{\shN}{\mathscr{N}}
\newcommand{\shL}{\mathscr{L}}
\newcommand{\shP}{\mathscr{P}}

\newcommand{\catC}{\mathcal{C}}

\newcommand{\aff}{\text{\rm aff}}
\newcommand{\Aff}{\text{\rm Aff}}
\newcommand{\alg}{\text{\rm alg}}

\newcommand{\Alb}{\operatorname{Alb}}

\newcommand{\Aut}{\operatorname{Aut}}

\newcommand{\can}{\text{\rm  can}}

\newcommand{\Cokernel}{\operatorname{Coker}}

\newcommand{\Ext}{\operatorname{Ext}}

\newcommand{\Frac}{\operatorname{Frac}}

\newcommand{\Gal}{\operatorname{Gal}}

\newcommand{\Hom}{\operatorname{Hom}}
\newcommand{\Hilb}{\operatorname{Hilb}}

\newcommand{\id}{{\operatorname{id}}}
\newcommand{\Image}{\operatorname{Im}}

\newcommand{\Kernel}{\operatorname{Ker}}

\newcommand{\invlim}{\varprojlim}

\newcommand{\lra}{\longrightarrow}

\newcommand{\maxid}{\mathfrak{m}}

\newcommand{\Num}{\operatorname{Num}}

 \newcommand{\NS}{\operatorname{NS}}
 
\newcommand{\primid}{\mathfrak{p}}
\renewcommand{\O}{\mathscr{O}}

\newcommand{\Pic}{\operatorname{Pic}}

\newcommand{\pr}{\operatorname{pr}}

\newcommand{\quadand}{\quad\text{and}\quad}

\newcommand{\ra}{\rightarrow}

\newcommand{\red}{{\operatorname{red}}}

\newcommand{\sep}{{\operatorname{sep}}}
\newcommand{\Set}{{\text{\rm Set}}}
\newcommand{\Sch}{\text{\rm Sch}}

\newcommand{\Spec}{\operatorname{Spec}}

\newcommand{\uHom}{\underline{\operatorname{Hom}}}

\begin{document}

\title[Para-abelian varieties and  Albanese maps]
      {Para-abelian varieties and  Albanese maps}

\author[Bruno Laurent]{Bruno Laurent}
\address{Heinrich Heine University D\"usseldorf, Faculty of Mathematics and Natural Sciences, Mathematical Institut,
40204 D\"usseldorf, Germany}
\curraddr{}
\email{bruno.laurent@hhu.de}

\author[Stefan Schr\"oer]{Stefan Schr\"oer}
\address{Heinrich Heine University D\"usseldorf, Faculty of Mathematics and Natural Sciences, Mathematical Institut,
40204 D\"usseldorf, Germany}
\curraddr{}
\email{schroeer@math.uni-duesseldorf.de}


\dedicatory{Revised version, 2 November 2023}

\begin{abstract}
We construct for every proper algebraic space over a ground field   an Albanese map to
a para-abelian variety, which is unique up to unique isomorphism. 
This holds   in the absence of rational points or ample sheaves,
and also for reducible or non-reduced spaces, under the mere assumption that the
structure morphism is in Stein factorization.
It also works under suitable assumptions in families. In fact  the treatment of the relative setting is crucial, 
even to understand the situation over ground fields. This also ensures that Albanese maps are equivariant with respect to actions of  group schemes.
Our approach depends on the notion of families of para-abelian varieties, where each geometric fiber   admits the  structure of an abelian variety, 
and   representability of   tau-parts in relative Picard groups, together with structure results on algebraic groups.
\end{abstract}
 
\maketitle
\tableofcontents

\newcommand{\dvs}{\text{\rm{dvs}}}
\newcommand{\prj}{\text{\rm{prj}}}
\newcommand{\afp}{\text{\rm{afp}}}
\newcommand{\trs}{\text{\rm{trs}}}
\newcommand{\abl}{\text{\rm{abl}}}

\newcommand{\rara}{\rightrightarrows}
\newcommand{\Rng}{\text{Rgn}}
\newcommand{\Opn}{\text{Opn}}
\newcommand{\catSO}{\operatorname{\mathcal{SO}}}
\newcommand{\Cls}{\text{Cls}}
\newcommand{\AffOpn}{{\rm\text{AffOpn}}}
\newcommand{\AlgSpc}{{\rm\text{\rm AlgSpc}}}
\newcommand{\catSF}{ \operatorname{\mathcal{SF}}}
\newcommand{\ab}{\text{\rm ab}}
\newcommand{\uExt}{\underline{\operatorname{Ext}}}

\section*{Introduction}
\mylabel{introduction}

The \emph{Albanese variety} and the \emph{Albanese map} are truly fundamental constructions in algebraic geometry.
The terms were   coined by Andr\'e Weil (\cite{Weil 1979}, page 438 and commentary on page 570, compare also \cite{Kleiman 2004}, page 428), in reference to 
work of Giacomo Albanese on zero-cycles and correspondences for  surfaces (\cite{Albanese 1934a} and \cite{Albanese 1934b}).
Originally, it was a purely \emph{transcendental construction}, depending on  path integrals over closed holomorphic one-forms.
A general formulation for compact complex spaces $X$ was given by Blanchard \cite{Blanchard 1956}:
The Albanese variety is the complex torus $V^*/\Delta$,  
constructed with the   dual   of the   vector space   $V\subset H^0(X,\Omega^1_{X/\CC})$ of closed forms, and
the Lie group closure $\Delta$   of the subgroup given by   integrals over loops. The resulting Albanese map 
$$
X\lra V^*/\Delta,\quad  x\longmapsto (\sigma\mapsto \int_{x_0}^x\sigma),
$$
is   defined by integrals over paths, and depends on the choice of a base point $x_0$.
For compact connected Riemann surfaces, this simplifies to $J=H^0(X,\Omega^1_X)^*/H_1(X,\ZZ)$.
Together with the principal polarization stemming from the intersection form this is called the \emph{jacobian variety}.
For genus $g\geq 1$ the Albanese map is  a closed embedding $X\subset J$, and for genus $g\geq 2$ one obtains an embedding
of the    Deligne--Mumford stacks $\shM_g\subset\shA_g$, which comprises   families of curves  and   principally polarized
abelian varieties, of   genus and dimension $g$, respectively.

It is natural to ask for  \emph{algebraic constructions} of   Albanese varieties and   Albanese maps, 
say for projective varieties, that extend  to general proper schemes $X$ over arbitrary ground fields $k$, and beyond. 
It quickly became clear that such a generalization is possible:
Matsusaka \cite{Matsusaka 1953} and Serre \cite{Serre 1960}   constructed the Albanese map  in   the classical language of algebraic varieties, by regarding
them as universal maps to abelian varieties, compare also the discussion of Esnault, Srinivas and Viehweg \cite{Esnault; Srinivas; Viehweg 1992}.
Grothendieck used  Picard schemes and Poincar\'e sheaves to obtain   Albanese maps, as outlined
in  \cite{FGA V}, and \cite{FGA VI}, Section 3.
The representability of the Picard functor was extended to arbitrary proper schemes by Murre \cite{Murre 1964}.
 
In the absence of rational points,  however, notorious complications arise  and  
Albanese maps take  values in   \emph{principal homogeneous spaces} rather than   abelian varieties.
These problems become even more pronounced over imperfect fields of characteristic $p>0$:
To our best knowledge, the existence of Albanese map, and its base-change behavior,
is only established for geometrically integral proper schemes, or more generally for geometrically reduced and 
geometrically connected proper schemes,
compare the discussions of Conrad \cite{Conrad 2017} and Wittenberg  \cite{Wittenberg 2008}.
Also see Achter, Casalaina-Martin and Vial \cite{Achter; Casalaina-Martin; Vial 2022} (formerly the appendix in the preprint \cite{Achter; Casalaina-Martin; Vial 2019}),
and \cite{Schroeer 2022} for non-proper $X$.

After the completion of the present work, we learned that Albanese maps where also constructed  in the setting of 
algebraic stacks by Brochard (\cite{Brochard 2021}, Section 7 and 8). They take values in commutative group stacks
that combine abelian varieties and finite group schemes, and exist under the condition that  $\Pic^\tau_{X/k}$ is proper
(loc.\ cit. Theorem 8.1).  The main advances of our work is the removal of the properness assumption, together with a  systematic 
development of the theory of para-abelian varieties and their relation to Picard scheme, which allows for a  
very geometric treatment of the subject.

The goal of this paper is to settle such issues, 
by systematically working in the \emph{relative setting} over some base scheme $S$,
and also to use  \emph{algebraic spaces} rather than schemes.
Recall that algebraic spaces are   important generalization of schemes  introduced by Michael Artin.
Roughly speaking,  they take over the role of   Moishezon spaces from complex geometry.  
However, their definition is much more indirect:
algebraic spaces are   contravariant functors $X:(\Aff/S)\ra(\Set)$ that satisfy the sheaf axiom
with respect to the \'etale topology, and are also otherwise  closely related to schemes.
Notions like the   underlying topological space
$|X|$ or the structure sheaf $\O_X$ still exist  but lose much of  their immediate significance, in comparison to the schematic situation.

Let $X$ be an algebraic space over $S$ 
such that the structure morphism  is proper, flat,  of finite presentation and  cohomologically flat
in degree $d=0$, with $h^0(\O_{X_s})=1$ for all points $s\in S$.  In particular, the structure morphism
$X\ra S$ is in Stein factorization. Our main result is:

\begin{maintheorem}
{\rm (see Thm.\ \ref{criterion albanese})} 
Assumptions as above. If   $\Pic^\tau_{X/S}$ admits a   family of maximal abelian subvarieties, then 
there is an Albanese map $f:X\ra\Alb_{X/S}$. Moreover, it is universal for morphisms into families of para-abelian varieties,
equivariant with respect to actions of group spaces, and commutes with base-change.
\end{maintheorem}

If the base scheme $S$ is integral, and the generic fiber of $\Pic^\tau_{X/S}\ra S$ is proper,
all assumptions are satisfied at least on some  dense open set $U$, over which the Albanese map
thus exists. If $S$ is   the spectrum of a field $k$ we get for each proper 
algebraic space $X$ with $h^0(\O_X)=1$ unconditional results. Note that $X$ may be reducible, non-reduced, or non-schematic, and the ground field $k$
is arbitrary. 
Building on this, the second author has treated the case of non-proper algebraic spaces over ground fields  
\cite{Schroeer 2022}.  Furthermore, the existence of  sign involutions on para-abelian varieties are studied in   \cite{Bergqvist; Dang; Schroeer 2023}.

To establish the above result we  develop  a theory of \emph{para-abelian varieties}, which is of independent interest.
Apparently, the name was coined by Grothendieck (\cite{FGA VI},  Theorem 3.3), but did not gain widespread use.
Again it is crucial to work in the relative setting: A \emph{family of para-abelian varieties}
is a    smooth proper morphism  $P\ra S$   such that for each point $s\in S$,
there is a field extension $k$ of the residue field $\kappa(s)$ such that $P \otimes k$
admits the structure of an abelian variety. \emph{Note that our definition does not involve any a priori torsor structure.}
However, we   show in Section \ref{Equivariance} that a certain inertia subsheaf in $\Aut_{P/S}$ is 
a family of abelian varieties, which a posteriori yields  a torsor structure.

The second main ingredient is a systematic study for the \emph{tau-part} $\Pic^\tau_{X/S}$ of the relative Picard functor.
Roughly speaking, it parameterizes invertible sheaves that are fiberwise numerically trivial.
We provide a new point of view, by directly using Artin's representability criteria \cite{Artin 1969}.
Here the representation is via algebraic spaces rather than schemes, even if $X$ is schematic, which highlights  that
the correct framework for our goals is given by algebraic spaces.
As already pointed out  by Grothendieck \cite{FGA VI}, the tau-parts are much better behaved than the connected components
$\Pic^0_{X_s/\kappa(s)}$  for the fiberwise Picard groups: The latter form an abelian sheaf that is usually not representable 
by algebraic spaces, and we give in Proposition \ref{enriques counterexample} an explicit example  with families of Enriques surfaces.

An important technical issue are the \emph{separation properties} of the ensuing quotient $N=\Num_{X/S}$ of the Picard space by its tau-part,
which we call the \emph{numerical sheaf}. These are algebraic spaces endowed with a group law whose fibers
are \'etale group schemes. In general, \emph{they are neither schematic nor separated, but   at least  locally separated},
which means that the diagonal monomorphism $N\ra N\times N$ is an embedding. This very weak separation axiom is enough to carry out
some crucial constructions that lead to Albanese maps.

An important insight of this paper is to \emph{re-define the notion of  Albanese maps} for    algebraic spaces $X$ over $S$:
Here it  denotes    pairs $(P,f)$  where $P$ is an $S$-family of para-abelian varieties, and $f:X\ra P$ is an $S$-morphism
such that for each point $s\in S$, the induced homomorphism  $f^*$ identifies  the abelian variety $\Pic^\tau_{P_s/\kappa(s)}$ with the 
\emph{maximal abelian subvariety} for the group scheme $\Pic^\tau_{X_s/\kappa(s)}$.
In some sense,  this definition is rather close to the original approach   with path integrals, and its  base-change properties follow  very easily.
We establish a posteriori that such Albanese maps have a universal property, which ensures uniqueness and equivariance.

In this approach, the \emph{maximal abelian subvariety} plays a crucial role. 
Extending some results of Brion \cite{Brion 2017} on the structure of group schemes $G$ of finite type over a base field $k$,
we introduce a \emph{three-step filtration} $G_3\subset G_2\subset G_1$ and describe the maximal abelian subvariety $G_\ab\subset G$
in terms of the  filtration and  extensions of abelian varieties by multiplicative groups, which appears to be of independent 
interest.

\medskip
The paper is organized as follows:
In Section \ref{Algebraic spaces} we briefly review algebraic spaces and local separatedness,
and give an existence result for quotients in a relative setting.
In Section \ref{Pic and tau} we thoroughly study  the tau-part $\Pic_{X/S}^\tau$ for the relative Picard functor, and
collect criteria for separatedness, properness,  flatness and smoothness.
This leads to the \emph{numerical sheaf} $\Num_{X/S}$, which is studied in Section \ref{Numerical sheaf}.
In Section  \ref{Para-abelian varieties} we develop a theory of para-abelian varieties and their families. 
Using actions of automorphism group schemes, we uncover in Section \ref{Equivariance} the relation between families of para-abelian varieties
and families of abelian varieties.
In Section \ref{Some extensions} and \ref{Maximal abelian} we work over ground fields $k$ and obtain   general results on the structure of group schemes
and the maximal abelian subvariety.
We then combine our results so far and state the definition of Albanese maps in Section \ref{Albanese maps}. Here we already 
establish some partial results concerning its uniqueness and existence. To proceed, we have to analyze the notion
of Poincar\'e sheaves in Section \ref{Poincare}.
In the final Section \ref{Existence}, we obtain existence and uniqueness for  Albanese maps and Albanese varieties
in full generality. In the Appendix, we collect some facts concerning embeddings of algebraic spaces.

\begin{acknowledgement}
We wish to thank Jeff Achter, Brian Conrad,  Laurent Moret-Bailly, Olivier Wittenberg and the referees for valuable comments
that helped to improve the paper and correct mistakes.
This research was conducted in the framework of the   research training group
\emph{GRK 2240: Algebro-geometric Methods in Algebra, Arithmetic and Topology}, which is funded
by the Deutsche Forschungsgemeinschaft. 
\end{acknowledgement}

\section{Algebraic spaces}
\mylabel{Algebraic spaces}
 
Throughout the   paper, $S$ usually denotes a base scheme. We write $(\Aff/S)$ for the category of affine  schemes $T=\Spec(R)$
endowed with  a structure morphism $T\ra S$.
Recall that an \emph{algebraic space} is a contravariant functor $X:(\Aff/S)\ra(\Set)$
satisfying the sheaf axiom with respect to the \'etale topology, such that   
the diagonal $X\ra X\times X$ is relatively representable by schemes,
and that there is an \'etale surjection $U\ra X$ from some scheme $U$. These are important
generalizations of schemes, because modifications, quotients, families, or moduli spaces  of schemes  
are frequently algebraic spaces rather than schemes.
We refer to the monographs of 
Olsson \cite{Olsson 2016}, Laumon and Moret-Bailly \cite{Laumon; Moret-Bailly 2000}, 
Artin \cite{Artin 1971}, Knutson \cite{Knutson 1971}, and to the stacks project
\cite{SP}, Part 0ELT. 

Note that there  is an  important novel separation axiom for
algebraic spaces, namely \emph{local separatedness}, which means that the diagonal $X\ra X\times X$
is an embedding (compare the Appendix for more details).
This is automatic for schemes, but  may fail for algebraic spaces (see the end of this section for examples).
Throughout we use the term \emph{embedding} in the sense of \cite{EGA I}, Definition 4.2.1, where the word \emph{immersion}
is used instead. Note that it is not necessarily closed, open, or quasicompact.
Also note that   some authors impose an additional condition, besides representability,  on      diagonals  for algebraic spaces, 
but we make no  separation assumptions whatsoever.

By definition, an algebraic space $X:(\Aff/S)\ra(\Set)$ satisfies the sheaf axiom with respect to the \'etale topology.
According to \cite{SP}, Lemma 076M it then also satisfies the sheaf axiom with respect to the fppf topology.
Throughout, we shall  regard $(\Aff/S)$ as a site with  this Grothendieck topology, if not stated otherwise.
It is obtained by the procedure explained in \cite{SGA 3a}, Expos\'e IV, Section 6 
from the families $(U_i\ra T)_{i\in I}$, where  $U_i\ra T$ are open embeddings, the index set $I$ is arbitrary, and $\bigcup U_i \ra T$ is surjective, 
together with the families $(V_j\ra T')_{j\in J}$, where  the index set $J$ is finite, the $V_j\ra T'$ are flat and of finite presentation,
and $\bigcup V_j\ra T'$ is surjective. Note that by dropping the assumption that $T$ and $U_i$ belong  to  $(\Aff/S)$,
the procedure allows to  pass to the larger category $(\Sch/S)$. 

Let $G$ be an algebraic space such that the $G(R)$, with $T=\Spec(R)$ from $(\Aff/S)$,  are endowed with   functorial group structures.
For simplicity, we say that $G$ is an \emph{algebraic space with  group structure}.
Let  $X$ be an algebraic space endowed with a $G$-action. One may form the quotient $X/G$ as a sheaf on $(\Aff/S)$.
Note that since the action is from the left, it would be more appropriate to write $G\backslash X$.
This common inconsistency should not cause confusion, because   
in most of our applications $G$ will be commutative.

We want to understand when the quotient $X/G$ is representable by an algebraic space.
Let us say that the  action $\mu:G\times X\ra X$ is  \emph{free} if  
$$
(\mu,\pr_2):G\times X\lra X\times X,\quad (\sigma,x)\longmapsto (\sigma x,x)
$$
is a \emph{monomorphism}. In other words, the groups $G(R)$ act freely on the sets $X(R)$, for all   $T=\Spec(R)$.
Note that this deviates from  other terminology that is sometimes used, where $(\mu,\pr_2)$ is assumed to 
be a \emph{closed embedding} (for example  \cite{Mumford; Fogarty; Kirwan 1993}, Definition 0.8).
For convenience of the reader, we state the following  result (\cite{SP}, Tag 06PH):
 
\begin{lemma}
\mylabel{quotient}
Suppose the   structure morphism  $G\ra S$ is  flat and locally of finite presentation,
and that the  action on $X$ is free.
Then the   quotient   $X/G$ in the category of sheaves is representable by an algebraic space. Moreover,
the quotient map $X\ra X/G$ is flat and locally of finite presentation, and the formation of $X/G$ commutes
with base-change. 
\end{lemma}

\medskip
In the situation of the Lemma, the quotient morphism $q:X\ra Y=X/G$ yields a cartesian diagram
$$
\begin{CD}
G\times X   @>(\mu,\pr_2)>>      X\times X\\
@Vq\circ\mu VV                        @VVq\times qV\\
Y           @>>\Delta>      Y\times Y.
\end{CD}
$$
The vertical map on the right is a torsor with respect to $(G\times G)_{Y\times Y}$, hence flat and locally of finite presentation.

\begin{proposition}
\mylabel{separation conditions}
Assumptions as in Lemma \ref{quotient}. Then  the monomorphism $(\mu,\pr_2):G\times X\ra X\times X$ is quasicompact, or an embedding,  or a closed embedding
if and only if the quotient $Y=X/G$ is quasiseparated, or locally separated,  or separated, respectively.
\end{proposition}

\proof
If $\Delta$ has one of the properties in question, so does the base-change $(\mu,\pr_2)$. 
For the converse, recall 
that for a morphism  being quasicompact, an   embedding,  or a closed  embedding is local in the range.
Choose an \'etale surjection $\amalg_{\lambda\in L} U_\lambda\ra X\times X$ for some affine schemes $U_\lambda$.
The compositions $U_\lambda\ra Y\times Y$ are flat and locally of finite presentation, hence universally open (\cite{EGA IVb}, Theorem 2.4.6).
Their images define open embeddings $V_\lambda\subset Y\times Y$ that cover $Y\times Y$, together with an fppf morphism $U_\lambda\ra V_\lambda$.
So if $(\mu,\pr_2)$ has one of the properties in questions, so does $\Delta$, by fppf descent.
\qed

\medskip
Let us close this section with two well-known examples of quotients   that are not locally separated, and in particular not schematic
(for related discussion see \cite{Olsson 2016},  Section 5.3 and \cite{SP}, Tag 02Z0 and \cite{Artin 1974}, Figure 1.1):

First suppose that $S=\Spec(k)$ is the spectrum of a field, that $H$ is a group scheme of finite type,
and $\Lambda\subset H(k)$ is   an infinite  subgroup. Regard the latter  as a constant group scheme $G=\amalg_{\sigma\in \Lambda} \Spec(k)$.
The canonical homomorphism $G\ra H$ is   a monomorphism, hence the translation action is free,
and the quotient $X=H/G$ exists as an algebraic space. It is not locally separated, because 
the monomorphism $G\ra H$ is not an embedding. 
Note that the  section $S\ra X$ coming from the neutral section in $H$ is not
an embedding, because this does not hold for the  base-change $G\ra H$, so  Corollary \ref{section embedding} below also ensures that 
$X$ is not locally separated. 
 
The second example lives over the spectrum $S=\Spec(R)$  of a discrete valuation ring. Write $U$ for the open set comprising the generic point,
and $Z$ for the complementary closed set.
Suppose $H$ is a relative  group scheme,
and $N$ is a closed subgroup scheme. Assume that the structure morphisms of $H\ra S$ and $N\ra S$ are fppf.   
Assume  that  $N_U\subset H_U$ is also open, and that $N_Z\subset H_Z$ is bijective.
Then the \emph{disjoint} union $G=N\amalg   (H_U\smallsetminus N_U)$ is a group scheme, and   the resulting    $G\ra H$ is a monomorphism. 
Obviously, the translation action of $G$ on $H$ is free,
and the structure morphism $G\ra S$ is fppf. It follows that the quotient $X=H/G$ exists as an algebraic space.
Now suppose there are $a,b\in H(S)$ such that $ab^{-1}$ does not lie in $G$, but $ab^{-1}|Z$ is  contained in $N(Z)$.
Then the restrictions of $a,b$ to both $Z$ and $U$ are congruent modulo  $G$. 
Now suppose that the diagonal $\Delta: X\ra X\times X$ would be an embedding. Then the cartesian diagram
$$
\begin{CD}
S'   @>>>    S\\
@VVV        @VV(\bar{a},\bar{b})V\\
X   @>>\Delta>    X\times X
\end{CD}
$$
defines a bijective embedding $S'\subset S$, whence $S'=S$.
This implies that $ab^{-1}\in G$, contradiction.
To make the situation concrete, one may choose $H=\mu_{p,R}$ and $N=\{e\}_R$ over the ring   $R=\ZZ_p[e^{2\pi i/p}]$.

\section{Pic and Pic-tau for families}
\mylabel{Pic and tau}
 
Let $S$ be a base scheme, $X$ be an algebraic space, and suppose that the structure morphism $ X\ra S$   is flat, proper, of finite presentation,
and cohomologically flat in degree $d=0$.
The latter ensures that the direct image of the structure sheaf $\O_X$ is locally free on $S$, and its formation  commutes with     base-change.
According to \cite{Artin 1969},  Theorem 7.3 the sheafification of $R\mapsto\Pic(X\otimes R)$
with respect to the fppf topology 
is representable by an algebraic space $\Pic_{X/S}$, which is locally of finite presentation,  quasiseparated and locally separated. 
By abuse of notation 
we here write $R\mapsto\Pic(X\otimes R)$ instead of the contravariant functor $T\mapsto \Pic(X\times_ST)$.
Note that in Artin's formulation of Theorem 7.3,   local separatedness is not stated explicitly because   only
algebraic spaces satisfying this separation axiom were considered.
Also note that $\Pic_{X/S}$ frequently fails to be schematic, and this already happens for certain
families of curves over discrete valuation rings (\cite{FGA VI}, Section 0).
For a nice overview in the realm of schemes, we also refer to  Kleiman's exposition \cite{Kleiman 2005}.

If $S $ is the spectrum of a field $k$, the component of the origin $\Pic_{X/k}^0$
is quasicompact, and the resulting quotient $\NS_{X/k}$ is a local system of finitely generated abelian groups,
called the \emph{N\'eron--Severi group scheme}.
Likewise, the preimage $\Pic_{X/k}^\tau$ of the torsion part in  $\NS_{X/k}$
is quasicompact. According to \cite{SGA 6}, Expos\'e XIII,  Theorem 4.6  the geometric points on $\Pic^\tau_{X/k}$
correspond to invertible sheaves that are numerically trivial, in other words the \emph{intersection number}
$(\shL\cdot C)=\chi(\shL)-\chi(\O_C) $
vanishes for every   curve $C\subset X$.
Note that the argument depends on Chow's Lemma, which holds in our context (\cite{Knutson 1971}, Theorem 3.1 or \cite{Rydh 2015}, Theorem 8.8), 
and the proof immediately carries over from schemes to algebraic spaces.

\emph{For a general base $S$, we now define $\Pic^\tau(X\otimes R)$ as the subgroup of the Picard group comprising
all invertible sheaves   that are fiberwise numerically trivial}. 
Clearly, this is functorial in $R$. We need the following result, which is already implicit in Artin's paper \cite{Artin 1969},
and was formulated by Brochard in the realm of stacks (\cite{Brochard 2012}, Theorem 3.3.3):

\begin{theorem}
\mylabel{pic tau}
The sheafification of $R\mapsto \Pic^\tau(X\otimes R)$ is representable by an algebraic space $\Pic^\tau_{X/S}$
whose structure morphism is of finite presentation (hence quasiseparated) and locally separated. 
Moreover, the inclusion into $\Pic_{X/S}$ is an open embedding.
\end{theorem}

\proof
Let us give a proof using Artin's representability criterion.
The problems are local in the base, so we may assume that $S=\Spec(A)$ is affine. 
Let $A_\lambda\subset A$ be the direct system of subrings that are finitely generated over the ring $\ZZ$.
Using the results from \cite{EGA IVc}, \S 8 there is an index $\mu$
such that $X=X_\mu\otimes_{A_\mu}A$ for some proper flat $X_\mu$ over $A_\mu$;
compare \cite{Rydh 2015} Proposition B.2 and B.3 for statements entirely in the realm of algebraic spaces.
We also have to ensure cohomological flatness:
Forming the \v{C}ech complex for some \'etale surjection $U_\mu\ra X_\mu$ with an  affine scheme $U_\mu$
and arguing as in \cite{Mumford 1970}, Section 5 we find a bounded complex $K^\bullet$ of
finitely generated projective $A_\mu$-modules giving an identification
$$
H^r(X_\mu\otimes B,\O_{X_\mu\otimes B}) = H^r(K^\bullet\otimes B),\quad r\geq 0
$$
that is functorial in the  $A_\mu$-algebras $B$. Clearly, the image  $M=\Image(K^0\ra K^1)$ is finitely generated
and its formation commutes with base-change. With the long exact Tor sequence for $0\ra H^0\ra K^0\ra M\ra 0$, we infer that the formation
of $H^0(X_\mu,\O_{X_\mu})$ commutes with base-change if and only if $M$ is locally free.
So by assumption, $M\otimes_{A_\mu} A$ is locally free. Hence there is some $\lambda\geq \mu$
such that $M\otimes_{A_\mu} A_\lambda$ is locally free.
Summing up, we may replace $X\ra \Spec(A)$ with the base change of $X_\mu\ra\Spec(A_\mu)$
to $A_\lambda$, and assume that $A$ is finitely generated over the excellent Dedekind ring $\ZZ$.

Artin established the representability for the sheafification of  $R\mapsto\Pic(X\otimes R)$ in \cite{Artin 1969}, Theorem 7.3  by an application of his  
Theorem 5.3, which involves checking certain conditions [0']--[5']. 
The reasoning for $\Pic^\tau(X\otimes R)$ is virtually the same. Only    condition [1'] requires   additional arguments:
Let $R$ be a complete local noetherian ring, with residue field $k=R/\maxid_R$.
We have to verify that the canonical map 
$$
\Pic^\tau(X\otimes R)\lra \invlim_n\Pic^\tau(X\otimes R/\maxid_R^{n+1})
$$
is bijective. It is injective, by Grothendieck's Existence Theorem 
(\cite{EGA IIIa}, Theorem 5.1.4 for schemes and \cite{SP}, 
Theorem 08BE for algebraic spaces).
Moreover, given invertible sheaves  $\shL_n$ on $X\otimes R/\maxid_R^{n+1}$ such
that the restriction of $\shL_{n+1}$ to $X\otimes R/\maxid_R^{n+1}$ is isomorphic to $\shL_n$,
and  that $(\shL_0\cdot C_0)=0$
for every integral curve $C_0$ in $X_0=X\otimes k$, the isomorphism classes come from some  invertible sheaf $\shL$ on $X$.
Our task is to check that the restriction to $X\otimes\kappa(\primid)$  is numerically trivial, for each prime ideal $\primid\subset R$.
By \cite{EGA II}, Proposition 7.1.4 it suffices to treat the case that $R$ is a discrete valuation ring,   $\primid$ is the zero ideal,
and that $S=\Spec(R)$. Let $\eta\in S$ be the generic point, and suppose that there is some integral curve $C_\eta\subset X_\eta$
with $(\shL_\eta\cdot C_\eta)\neq 0$.
The closure $C\subset X$ of this curve is flat, whence  the closed fiber 
$C_0\subset X_0$ is a curve. Note that by  \cite{EGA IVd}, Proposition 21.9.11 the total space $C$ carries an ample sheaf, 
hence the algebraic space $C$ is  a scheme. Since Euler characteristics are constant in families (\cite{EGA IIIb}, Theorem 7.9.4), we see 
 $(\shL_0\cdot C_0)=(\shL_\eta\cdot C_\eta)\neq 0$, contradiction.

This shows that $\Pic^\tau_{X/S}$ is an algebraic space.
According to \cite{SGA 6}, Expos\'e XIII, Theorem 4.7 the monomorphism to $\Pic_{X/S}$ is an open embedding,
and the structure morphism $\Pic^\tau_{X/S}\ra S$ is quasicompact. The proof relies on   approximation and Chow's Lemma as above,
and carries over from schemes to algebraic spaces. Since $\Pic_{X/S}$ is   locally of finite presentation, quasiseparated and locally separated,
the same holds for the open subspace $\Pic^\tau_{X/S}$.
\qed

\medskip
In \cite{FGA VI}, Grothendieck defined   $\Pic^\tau_{X/S}$ as the   subsheaf comprising those $R$-valued points
of $\Pic_{X/S}$ that are torsion in the fiberwise N\'eron--Severi groups, and showed in his Theorem 1.1
that it is open, under the assumption that $\Pic_{X/S}$ is a scheme.
The above approach seems to be more adequate in the realm of  algebraic spaces.

Similarly, one defines $\Pic^0_{X/S}$ as the \emph{abelian subsheaf} comprising the $R$-valued points of $\Pic_{X/S}$ that
are trivial in the fiberwise N\'eron--Severi groups. This indeed seems the only possible approach,
because the subsheaf is not representable in general, as the following counterexample shows:

Suppose $S$ is the spectrum of $\FF_2[[t]]$, 
and let $X\ra S$ be a family of Enriques surfaces whose generic fiber is classical but whose
closed fiber is ordinary or supersingular (for details see \cite{Bombieri; Mumford 1976}, Section 3). 
Then $G=\Pic_{X/S}^\tau$ is finite and flat of degree $n=2$. It is the union  $G=G'\cup G''$
of two sections  intersecting in the closed fiber, which is a singleton. 

\begin{proposition}
\mylabel{enriques counterexample}
In the above situation, the  subsheaf $\Pic^0_{X/S}$ inside $\Pic^\tau_{X/S}$ is not representable by an algebraic space.
\end{proposition}

\proof 
Suppose $H=\Pic^0_{X/S}$ is representable. By definition, the canonical morphism to $G=\Pic^\tau_{X/S}$
is a monomorphism, and the set of image points is the zero-section.
It follows that the structure morphism $H\ra S$ is a universal homeomorphism.
In particular, it is universally closed, separated, and has affine fibers.
According to \cite{Rydh 2015}, Theorem 8.5 it must be integral, and in particular affine. Write $H=\Spec(A)$
and $G=\Spec(R)$, and consider the canonical homomorphism $R\ra A$. 
The induced map $R/t^nR\ra A/t^nA$ is bijective for each $n\geq 0$. We have $\bigcap_{n\geq 0} t^nA=0$ by Krull's Intersection Theorem, thus $R\ra A $ is injective.
It follows that the image of $H\ra G$ contains both generic points, contradiction.
\qed

\medskip
Let us finally collect the basic properties of the tau-part:

\begin{proposition}
\mylabel{properties pic tau}
The algebraic space $\Pic^\tau_{X/S}$ has the following properties:
\begin{enumerate}
\item
If the geometric fibers of $f:X\ra S$ are   integral, then  the structure morphism $\Pic^\tau_{X/S}\ra S$ is   separated.
\item
If the   geometric fibers of $f:X\ra S$ are integral and   locally factorial, then  
$\Pic^\tau_{X/S}\ra S$ is proper.
\item
If for all  points $s\in S$ we have $h^1(\O_{X_s})- h^2(\O_{X_s})=b_1(X_s)/2$ then the morphism $\Pic^\tau_{X/S}\ra S$ is flat.
\item
Suppose for each Artin local ring $B$ with residue field $k=B/\maxid_B$, and each   $\Spec(B)\ra S$,
there is an algebraic space $Z$, some proper $B$-morphisms $h_1,\ldots,h_r:Z\ra X\otimes B$ and integers $n_1,\ldots,n_r$ such that the induced map
\[
\sum n_ih_i^*:H^2(X\otimes k ,\O_{X\otimes k})\lra H^2(Z\otimes k,\O_{Z\otimes k})
\]
is injective, whereas
$\sum n_ih_i^*:\Pic^\tau_{X\otimes B/B}\ra\Pic^\tau_{Z/B}$ is zero.
Then $\Pic^\tau_{X/S}\ra S$ is smooth.
\end{enumerate}
\end{proposition}
 
\proof
It suffices to treat the case that $S=\Spec(R)$ is the spectrum of a noetherian ring.
Set $G=\Pic^\tau_{X/S}$.
We start with assertion  (i) and (ii), which follow from   Valuative Criteria, similar to 
\cite{FGA VI}, Theorem 2.1. 
Suppose the fibers are geometrically integral. 
Then the map $\O_S\ra f_*(\O_X)$ is bijective.
Let $V$ be a discrete valuation ring
with field of fractions $F=\Frac(V)$, such that the residue field $k=V/\maxid_V$ is algebraically closed
and that $V$ is complete.  

For (i) we have to verify 
that the map $G(V)\ra G(F)$ is injective (\cite{Laumon; Moret-Bailly 2000}, Proposition 7.8 or \cite{SP}, Lemma 03KV).
Let $l_1,l_2\in G(V)$ be two elements that coincide in $G(F)$.
Without restriction, we may assume $V=R$. Then $H^2(S,\GG_m)=0$, for example by \cite{Milne 1980}, Chapter IV, Corollary 1.7 combined with Corollary 2.12. 
Now the Leray--Serre spectral sequence for the structure morphism $X\ra S$   shows that 
the canonical map $\Pic^\tau(X)\ra G(R)$ is bijective and that the canonical map $\Pic^\tau(X\otimes F)\ra G(F)$ is injective.
So the sections $l_i$ come from invertible sheaves
$\shL_i$ that are isomorphic on the complement $U$ of the integral Cartier divisor $D=X\otimes k$.
Choose an isomorphism $\varphi:\shL_1|U\ra \shL_2|U$.
It extends to  a homomorphism $\varphi:\shL_1(-nD)\ra\shL_2$
for some integer $n\geq 0$, as in \cite{EGA I}, Theorem 6.8.1. This gives a short exact sequence
$0\ra \shL_1(-nD)\ra \shL_2\ra \shF\ra 0$.
The cokernel  is an invertible sheaf on some effective Cartier divisor $E\subset X$
supported by $D$. Since the latter is integral, we have $E=mD$ for some integer $m\geq 0$.
Setting $\shN=\shL_2\otimes\shL_1^{\otimes-1}\otimes\O_X(rD)$ with $r=m-n$,  we get another
short exact sequence $0\ra \shN(-mD)\ra \shN\ra \shF\ra 0$.
Computing the determinant of $\shF$  with the second sequence  and using $\O_X(D)\simeq\O_X$ gives
$\det(\shF)\simeq\O_X$. Using the first sequence then reveals  $\shL_1\simeq\shL_2$. This shows (i). Note that this generalizes a result
of Ischebeck \cite{Ischebeck 1979}, who considered affine schemes.

For assertion (ii), suppose we have a point $l\in G(F)$. Choose a finite separable extension $F\subset F'$ so 
that the image $l'\in G(F')$ comes from an invertible sheaf   on $X\otimes F'$, and let
$V'$ be the integral closure of $V$. Our task is to check that $l'$ lies in the subset $G(V')\subset G(F')$
(\cite{Laumon; Moret-Bailly 2000}, Theorem 7.3 or \cite{SP}, Lemma 0A3Z).
Without loss of generality, we may assume that $R=V=V'$. 
Let $\shL_F$ be the invertible sheaf on $X\otimes F$ corresponding to the point $l$,
and $\shF$ be any coherent extension to $X$. We now check that its bi-dual $\shL$
is invertible.
The total space $X$ is integral, because this holds for the closed fiber $X\otimes k$.
Since $X\otimes F$ is locally factorial, the same holds for $X$, by Nagata's result (\cite{Nagata 1957}, Lemma 1).
In particular, $X$ is normal, and $\shL$ is reflexive of rank one, hence invertible.
This settles (ii). 

Assertion (iii) is due to Ekedahl, Hyland and Shepherd-Barron (\cite{Ekedahl; Hyland; Shepherd-Barron 2012}, Proposition 4.2).
Note the   Betti numbers $b_i(X_s)$ are  vector space dimensions for
the \'etale cohomology groups 
$H^i(X_{\bar{s}},\QQ_l(i)) = \invlim_{\nu} H^i(X_{\bar{s}},\mu_{l^\nu}^{\otimes i})\otimes_{\ZZ_l}\QQ_l
$,
which by definition are computed over the algebraic closure of the residue field $\kappa(s)$.

We finally  come to statement (iv), which is an abstraction of Mumford's arguments for abelian varieties
(\cite{Mumford; Fogarty; Kirwan 1993}, Proposition 6.7).
 We saw in Theorem \ref{pic tau} that $G=\Pic_{X/S}^\tau$ is of finite presentation
over $S$. By   definition of smoothness, we have to check that the canonical map
$G(A)\ra G(A/J)$ is surjective, for each $R$-algebra $A$ and each ideal $J$ with $J^2=0$.
In light of \cite{EGA IVd}, Remark 17.5.4 it suffices
to consider the case where $A$ is a local Artin ring. We may further assume that 
the residue field $A/\maxid_A$ is algebraically closed,   that the ideal $J$ has length one,
and that $A=R$. 

Let $h:X\ra S$ be the structure morphisms and $X^\aff=\Spec \Gamma(X,\O_X)$.
The proof for \cite{Schroeer 2020}, Lemma 1.4 reveals that we have an identification
$H^2(S, h_*(\GG_{m,X})) = H^2(X^\aff,\GG_m)$. The Brauer group $H^2(X^\aff,\GG_m)$ vanishes,   because  $X^\aff$ is
an Artin scheme with algebraically closed residue fields. Consequently, the Leray--Serre spectral sequence
for $h$ gives an exact sequence
$$
\Pic^\tau(X)\lra G(R)\lra H^2(S, h_*(\GG_{m,X}))=0.
$$
We infer that  the canonical map  $\Pic^\tau(X)\ra G(R)$ and   also $\Pic^\tau(X\otimes B)\ra G(B)$ is bijective, with $B=A/J$.
The short exact sequence $0\ra \O_{X\otimes k}\ra \O_X^\times\ra \O_{X\otimes B}^\times \ra 1$
and the corresponding sequence for $Z$ induce a commutative diagram
\[
\begin{CD}
\Pic^\tau(X) @>>>       \Pic^\tau(X\otimes B)    @>>>       H^2(X\otimes k,\O_{X\otimes k})\\
@VVV                    @VVV                                @VVV            \\
\Pic^\tau(Z) @>>>       \Pic^\tau(Z\otimes B)    @>>>       H^2(Z\otimes k,\O_{Z\otimes k}).
\end{CD}
\]
Here the vertical arrows are given by the linear combination $\sum n_ih_i^*$.
By assumption, the vertical map on the right is injective, but the vertical map in the middle is zero.
So by  exactness of the upper row, each invertible sheaf on $X\otimes B$ that is fiberwise numerically trivial extends to $X$.
\qed

\medskip
Note that assertion (i) actually hold for the structure morphisms $\Pic_{X/S}\ra S$ of the whole Picard scheme,
with the same proof.

\section{The numerical sheaf}
\mylabel{Numerical sheaf}

Let $S$ be a base scheme, and $X$ be an algebraic space whose structure morphism $X\ra S$ is proper,
flat, of finite presentation and cohomologically flat in degree $d=0$.
We then define an abelian sheaf $\Num_{X/S}$ on the site $(\Aff/S)$ by the short exact sequence
$$
0\lra \Pic^\tau_{X/S}\lra \Pic_{X/S}\lra \Num_{X/S}\lra 0,
$$
and call it the \emph{numerical sheaf}.
Note that the translation action of $\Pic^\tau_{X/S}$ on $\Pic_{X/S}$ is free, so the formation of the quotient
commutes with base-change. In particular, for each point $s$ and each field $k'$ containing the separable closure $\kappa(s)^\sep$,
the  abelian group $\Num_{X/S}(k')$
is free, and its rank is the Picard number $\rho\geq 0$ of $X\otimes \kappa(s)^\alg$, by \cite{SGA 6}, Expos\'e XIII, Theorem 5.1.
We actually have:

\begin{theorem}
\mylabel{num representable}
Suppose $\Pic^\tau_{X/S}\ra S$ is flat. Then  $\Num_{X/S}$ is representable by an algebraic space. The structure morphism $\Num_{X/S}\ra S$
is locally of finite presentation,  quasiseparated and locally separated. Moreover, all  fibers are separated, schematic and \'etale.
\end{theorem}

\proof
The problem is local in $S$, and we may assume that $S$ is the spectrum of a noetherian ring $R$.
The first assertion follows from Lemma \ref{quotient}. Since $\Pic^\tau_{X/S}$ is noetherian,
the open embedding into $\Pic_{X/S}$ is  quasicompact, and the   assertion on the structure morphism follows from 
Proposition \ref{separation conditions}.
To establish the statement on the fibers, we may assume that $R=k$ is a field.
Then $G=\Pic_{X/k}$ is a scheme, and $H=\Pic^\tau_{X/k}$ is an open subgroup scheme. Hence it is also closed,
and the zero element $e\in G/H$ is a closed point. 
In particular, the algebraic space with group structure $G/H$ is separated. 
It must be \'etale and hence zero-dimensional, because $H\subset G$ is open.
Summing up, the algebraic space $G/H$ is locally of finite type over the field $k$, and zero-dimensional.
It the must be schematic by  \cite{SP}, Lemma 06LZ.
\qed

\medskip
%
Suppose now  that $\Pic^\tau_{X/S}$ is flat over $S$, so that the numerical sheaf $\Num_{X/S}$ is an algebraic space,
for example by the criteria  in Proposition \ref{properties pic tau}.
We need the following fact, which is a variant  of Mumford's Rigidity Lemma (\cite{Mumford; Fogarty; Kirwan 1993}, Proposition 6.1):

\begin{lemma}
\mylabel{rigidity}
Let $A$ be an algebraic space whose structure morphism $\alpha:A\ra S$ is proper, flat, of finite presentation,
cohomologically flat in degree $d=0$, and with $h^0(\O_{A_s})=1$ for all $s\in S$.
Then for each   $\varphi\in\Num_{X/S}(A)$, there is a unique $\sigma\in \Num_{X/S}(S)$ making the following diagram commutative:
\begin{equation}
\label{rigidity diagram}
\begin{tikzcd} 
A\ar[dr,"\alpha"']\ar[rr,"\varphi"]      &       & \Num_{X/S}\ar[dl]\\
        & S\ar[ur, bend right=30, "\sigma"']\\ 
\end{tikzcd}
\end{equation} 
\end{lemma}

\proof
Write $N=\Num_{X/S}$.
In light of the uniqueness assertion, the problem is local, and as usual we   may assume that $S$ is the spectrum
of a noetherian ring $R$.  Again by uniqueness, together with fppf descent 
(\cite{SGA 1}, Expos\'e VIII, Theorem 5.2 for schemes and \cite{SP}, 
Lemma 0ADV for algebraic spaces),
it suffices to treat the case $\alpha:A\ra S$ admits a section $\tilde{\sigma}$. 
So  the diagram \eqref{rigidity diagram} reveals that if $\sigma$ exists, it must coincide with the section $\varphi\circ\tilde{\sigma}$,
which thus already settles the uniqueness assertion. It remains to verify that $\sigma=\varphi\circ\tilde{\sigma}$
indeed makes the diagram commutative.

To start with, suppose that $R=k$ is a field, so that $N$ is a separated scheme whose underlying topological space is discrete.  Since $A$ is proper,
the set-theoretical image $\varphi(A)\subset N$ is closed and quasicompact. It carries exactly one scheme structure,
because $N$ is \'etale, and must be  the spectrum of an \'etale $k$-algebra $L=k_1\times\ldots\times k_r$.
In particular, the schematic image is affine. Using $h^0(\O_A)=1$ we conclude $L=k$, and $\varphi=\sigma\circ \alpha $ follows.

For the general case, consider the algebraic space $U$ defined by the cartesian diagram
$$
\begin{CD}
U	@>>>				N\\
@VVV					@VV\Delta V\\
S	@>>(\sigma,\varphi\circ\tilde{\sigma})>	N\times N.
\end{CD}
$$
It follows from \cite{EGA IVd}, Corollary 17.4.2, together with Theorem \ref{num representable} that the diagonal $\Delta:N\ra N\times N$ is an open embedding.
In turn, the same holds for $U\ra S$. By construction, this open set comprises all points at which $\sigma$ and $\varphi\circ\tilde{\sigma}$ coincide.
With the preceding paragraph we infer $U=S$, and thus $\sigma=\varphi\circ\tilde{\sigma}$.
\qed

\section{Para-abelian varieties}
\mylabel{Para-abelian varieties}

Let us call an  algebraic space $P$ over some field $k$ a \emph{para-abelian variety} if there
is a field extension $k\subset k'$ such that the base-change $P'=P\otimes_kk'$ admits
the structure of an abelian variety.   
The terminology goes back to Grothendieck, who  introduced  
it in \cite{FGA VI},  Theorem 3.3  by a different condition. A posteriori, we shall see  that
the notions are equivalent. 
By fpqc descent,  our $P$ is proper and smooth over $k$, with $h^0(\O_P)=1$. 
Moreover:

\begin{lemma}
\mylabel{projective}
The algebraic space $P$ is a projective scheme.
\end{lemma}

\proof
Choose an algebraically closed field extension $k'\subset\Omega$.
According to \cite{Mumford 1970}, page 62 there is an ample sheaf  $\shL$  on the base-change $P\otimes_k\Omega$. 
This defines a morphism $\Spec(\Omega)\ra\Pic_{P/k}$,
which factors over some connected component   $\Pic_{P/k}^{l}$.
Fix a closed point $a$ in this component. Then the field extension $k\subset\kappa(a)$ is finite,
and there is a finite extension $\kappa(a)\subset k''$ and some invertible sheaf $\shA$ on $P\otimes k''$
mapping to $a\in \Pic_{P/k}$. Choose an embedding $k''\subset \Omega$. 
By the numerical criterion for ampleness the base-change  of $\shA$ to $P\otimes\Omega$ is  ample  (\cite{Kleiman 1966}, page 343, Theorem 1).
It follows that   $P\otimes k''$ admits an ample invertible sheaf $\shN''$. Let $\shN=N (\shN'')$ be the norm
with respect to the finite locally free morphism $P\otimes k''\ra P$.
According to \cite{EGA II}, Proposition 6.6.1 the base-change  of $\shN$ to $P\otimes k''$ is ample.
Thus $\shN$ is ample on $P$.
\qed

\medskip
For our purposes it is crucial to work in the relative setting: 

\begin{definition}
\mylabel{def family para-abelian}
A \emph{family of para-abelian varieties} over some scheme $S$
is an algebraic space $P$, together with a morphism $P\ra S$ that is proper, flat and of finite presentation, such that the fibers $P_s$ are para-abelian varieties
over the residue field $\kappa(s)$, for every $s\in S$.
\end{definition}

Particular examples are the   \emph{families of abelian varieties}. By this we mean an algebraic space $A$,
together with a proper flat  morphism of finite presentation $A\ra S$ endowed with a group structure, such that all fibers are abelian varieties.
These are often called \emph{abelian schemes} in the literature.
According to  Raynaud's result (see \cite{Faltings; Chai 1990}, Theorem 1.9)
the total space $A$ indeed must be a scheme. It actually satisfies the \emph{AF-property}, that is, each finite set of points
admits a common affine open neighborhood, provided the base scheme $S$ is affine. 
Actually, the total space is quasiprojective if the base is affine and normal (\cite{Faltings; Chai 1990}, Remark 1.10).
Note, however, that there are examples without ample sheaves
(\cite{Raynaud 1970}, Chapter XII, 4.2), and this happens already over the spectrum of the   ring of dual numbers $R=\CC[\epsilon]$.

In what follows, we fix a family of para-abelian varieties $f:P\ra S$. Note that the  structure morphism
is    surjective and smooth.
Moreover, the canonical map $\O_S\ra f_*(\O_P)$ is bijective (\cite{EGA IIIb}, Proposition 7.8.6).
In other words, $f$ is cohomologically flat in degree $d=0$, and $h^0(\O_{P_s})=1$ for
all points $s\in S$.
In contrast to families of abelian varieties, the total space $P$ often fails to be a scheme.
This already happens in relative dimension $g=1$ over local schemes $S$ of dimension $n=2$, see \cite{Raynaud 1970}, Chapter XIII, Section 3.2
and also \cite{Zomervrucht 2015}.

We seek to relate the sheaves  $\Aut_{P/S}$ and   $\Pic_{P/S}$   to $P\ra S$. Let us start with some useful observations,
which generalize  \cite{Mumford; Fogarty; Kirwan 1993}, Theorem 6.14:

\begin{proposition}
\mylabel{group law}
For each   $e\in P(S)$, there is a unique group law $\mu:P\times_SP\ra P$
that turns $P\ra S$ into a family of abelian varieties, with $e:S\ra P$ as the zero section.
\end{proposition}

\proof
Uniqueness and hence also existence are local problems, so we may assume that $S=\Spec(R)$ is affine.
Suppose there are two group laws $\mu_1$ and $\mu_2$ with $e$ as zero section.
Recall that the algebraic space $P$ then must be schematic.
Since $P$ and hence $P\times_SP$ are of finite presentation,  the scheme $P$ and the morphisms $\mu_i$
are already defined over some noetherian subring $R_0\subset R$.
Now \cite{Mumford; Fogarty; Kirwan 1993}, Corollary 6.6 ensures that $\mu_1=\mu_2$. This settles uniqueness. 

It remains to verify existence. For this it also suffices to treat the case that $R$ is noetherian.
Suppose first that there is an fpqc extension $R\subset R'$ such that a group law $\mu'$ 
exists for $P'=P\otimes_RR'$, with some origin $e'\in P(R')$. Using translation by $e\otimes 1-e'$,
we may assume that $e'$ is the base-change of $e$. Consider the ring $R''=R'\otimes_RR'$.
By fpqc descent, we have to verify that the two pull-backs $\mu'\otimes 1$ and $1\otimes \mu'$ to $R''$ coincide.
Both are group laws, and in both cases the origin is the pull-back of $e$.
Uniqueness ensures $\mu'\otimes 1 = 1\otimes\mu'$.
Note that this settles the assertion  if $R=k$ is a field.

Suppose next that $S$ is the spectrum of some  local Artin ring $R$.
Then $P_\red $ is the closed fiber, which is schematic according to  Proposition \ref{projective}.
By \cite{Rydh 2015}, Corollary 8.2 the total space $P$ is schematic as well. (This result already
appears in \cite{Knutson 1971}, Theorem 3.3, at least for quasiseparated algebraic spaces.)
Hence the group law $\mu$ exists by  \cite{Mumford; Fogarty; Kirwan 1993}, Proposition 6.15.

Now suppose that $R$ is a general noetherian ring.
Fix a closed  point $a\in S$, corresponding to a maximal ideal $\maxid\subset R$.
The preceding paragraph gives a formal group law over the completion $\hat{R}=\invlim R/\maxid_a^n$.
It comes from a  group law over   $\hat{R}$, according to 
Grothendieck's Existence Theorem (\cite{EGA IIIa}, Theorem  5.4.1 for schemes
and \cite{SP}, 
Lemma 0A4Z for algebraic spaces).
Since $R_\maxid\subset \hat{R}$ is an fpqc extension, the group law already exists over $R_\maxid$.
Thus it is already defined over some open neighborhood $U\subset S$ of $a$
(\cite{EGA IVc}, Theorem 8.8.2 for schemes and \cite{Rydh 2015}, Proposition B.2 for algebraic spaces).
Applying this for all closed points $a\in S$, we obtain an open covering $S=U_1\cup\ldots\cup U_r$
such that the group law exists over each $U_i$. By uniqueness, these  local group laws glue and yield
the desired global group law.
\qed
 
\begin{corollary}
\mylabel{stricly henselian}
Suppose $S=\Spec(R)$ is   henselian, with closed point $a\in S$. Then $P\ra S$ admits the structure of a  family of abelian varieties if and only if 
the closed fiber $P_a$ contains a rational point. 
\end{corollary}

\proof
The condition is obviously necessary. Conversely, suppose there is a rational point $e_a\in P_a$.
It follows from \cite{EGA IVd}, Corollary 17.16.3 that 
there is a subscheme $Z\subset P$ containing $e_a$, and such that $Z\ra S$ is \'etale and quasi-finite.
In turn, the singleton $\{e_a\}$ is a connected component of the closed fiber $Z_a$. 
It corresponds to a connected component $U\subset Z$, because $R$ is henselian. Moreover, $U\ra S$ is \'etale and finite,
and thus defines a section $e\in P(S)$ extending $e_a$. By the Proposition, $P\ra S$ becomes a family of abelian varieties.
\qed

\begin{corollary}
\mylabel{etale surjection}
There is an \'etale surjection $S'\ra S$ such that the base-change $P'=P\times_SS'$ admits the structure of 
a family of abelian varieties over $S'$.
\end{corollary}

\proof
According to \cite{EGA IVd}, Corollary 17.16.3 there is an \'etale surjection $S'\ra S$ such that $P(S')$ is non-empty.
The Proposition ensures the existence of a group law.
\qed

\medskip
Recall that $B=\Pic^\tau_{P/S}$ is an algebraic space, and the morphism $B\ra S$ is of finite presentation. For each integer $n\geq 1$,
the kernel $B[n]$ for multiplication by $n$ is another algebraic space. 
Write $G_n=\Aut_{B[n]/S}$ for the ensuing automorphism sheaf.

\begin{corollary}
\mylabel{relative affine}
In the above situation, the structure morphism $B\ra S$ is a family of abelian varieties. 
Moreover,  for each integer $n\geq 1$ the   morphism  $B[n]\ra S$ is finite and  locally free,
and   $G_n\ra S$ is relatively representable by affine schemes. 
\end{corollary}

\proof
We first check that $B\ra S$ is a family of abelian varieties.
The  morphism is proper  according to Proposition \ref{properties pic tau}, part (ii).
To check smoothness, it suffices to treat the case that  $P\ra S$ admits a section (\cite{EGA IVd}, Proposition 17.7.1).
So $P\ra S$ becomes a family of abelian varieties, by our Proposition.
Mumford then showed   that $B=\Pic^\tau_{P/S}$ is smooth (\cite{Mumford; Fogarty; Kirwan 1993}, Proposition 6.7).
He actually assumed that $P\ra S$ is projective to have the existence of $\Pic^\tau_{B/S}$ as   projective scheme.
However, his arguments   carry over without change to our situation (relying on Theorem \ref{pic tau} and Proposition \ref{properties pic tau}, part (iv)).
It remains to verify that $B\ra S$ has geometrically connected fibers.
For this it suffices to treat the case that $S$ is the spectrum of an algebraically closed field $k$, and
that $P$ is an abelian variety.  Then the N\'eron--Severi group
$\NS(P)$ is torsion free, according to \cite{Mumford 1970}, Corollary 2 on page 178.
Thus $B\ra S$ is a family of abelian varieties.

For the remaining assertions, it suffices to treat the case that $S$ is the spectrum of a noetherian ring $R$,
and $P\ra S$ is a family of abelian varieties. Since $B\ra S$ is separated, the inclusion $B[n]\subset B$ is closed,
hence the structure morphism $B[n]\ra S$ is proper. It is also quasi-finite (\cite{Mumford 1970},
Appendix to \S 6).
By \cite{EGA IVd}, Corollary 18.12.4   it must be finite.
For each point $s\in S$, the fiberwise multiplication $n:B_s\ra B_s$ is finite and surjective,
and hence flat (\cite{Serre 1965}, Proposition 22 on page IV-37).
In light of \cite{EGA IVc}, Proposition 11.3.11 this holds true for $n:B\ra B$.
Consequently $B[n]\ra S$ is finite and locally free.
It follows that $G_n\ra S$ is affine, for example by \cite{Lorenzini; Schroeer 2020}, Lemma 4.1.
\qed

\medskip
Let $G $ be an algebraic space endowed with a group structure,  and assume that the structure morphism $g:G\ra S$ is   flat and of finite presentation,
and that $\O_S\ra g_*(\O_G)$ is bijective.
Suppose that we have a relative $G$-action  on $P$. By functoriality, it induces a relative action on $B=\Pic^\tau_{P/S}$.

\begin{corollary}
\mylabel{induced action trivial}
In the above situation, the $G$-action on $B$ is trivial. 
\end{corollary}

\proof
It suffices to treat the case that $S$ is the spectrum of a noetherian ring $R$.
The induced action on $B[n]$ is trivial, because the latter are affine, whereas $R=\Gamma(G,\O_G)$.
According to \cite{EGA IVc}, Theorem 11.10.9  the collection of closed subgroup schemes $B[n]\subset B$, $n\geq 1$ is schematically dense, 
and it follows by from loc.\ cit., Proposition 11.10.1 that the action on $B$ must be trivial as well.
\qed

\section{Equivariance}
\mylabel{Equivariance}

Fix a base scheme $S$, and let $X$ be an algebraic space whose structure morphism  $X\ra S$   is  locally of finite presentation and separated.
According to \cite{Artin 1969}, Theorem 6.1 the Hilbert functor $\Hilb_{X/S}$ is representable
by an algebraic space that is locally of finite presentation and separated. Recall that its $R$-valued points
are the closed subspaces $Z\subset X\otimes R$ such that the projection $Z\ra \Spec(R)$ is proper, flat and
of finite presentation.
If furthermore $X\ra S$ itself is proper and flat, one sees that $\Aut_{X/S}$ is an open subspace of $\Hilb_{(X\times X)/S}$,
by interpreting automorphisms via their graphs.

Now let $P\ra S$ be a family of para-abelian varieties. 
The action from the left of the algebraic  space 
$\Aut_{P/S}$ on $P$ induces an action from the right on the family of abelian varieties $B=\Pic^\tau_{P/S}$, via pull-back of invertible sheaves.
Write $G\subset\Aut_{P/S}$ for the ensuing \emph{inertia subgroup sheaf}; its group of $R$-valued points comprises the $R$-isomorphisms  $f:P\otimes R\ra P\otimes R$
where the induced map $f^*:B\otimes R\ra B\otimes R$ is the identity.

\begin{proposition}
\label{automorphism group}
\mylabel{inertia abelian}
The  inclusion $G\subset\Aut_{P/S}$ is representable by   open-and-closed embeddings, the structure morphism
$G\ra S$ is a family of abelian varieties, and the total space $G$ is a scheme.
\end{proposition}

\proof
First observe that once we know that $G\ra S$ is a family of abelian varieties, the total space must be a scheme by
Raynaud's result (\cite{Faltings; Chai 1990}, Theorem 1.9).
To verify the statements on $G\subset \Aut_{P/S}$ and $G\ra S$, it suffices to treat  the case that $P\ra S$ is a family of abelian varieties, 
and $S$ is the spectrum of a noetherian ring $R$.
Note that then  the translation action gives a monomorphism $P\ra\Aut_{P/S}$.
This is a closed embedding, because $P$ is proper and $\Aut_{P/S}$ is separated
(\cite{EGA IVd}, Corollary 18.12.6).

By Corollary \ref{induced action trivial}  we have $P\subset G$ as subsheaves inside  $\Aut_{P/S}$.
We claim that this inclusion is an equality.
This is a statement on $R'$-valued points; by making 
a base-change it suffices to treat the case $R=R'$.
Suppose we have some automorphism $f:P\ra P$ that lies in $G(R)$.
Write $f=t_a\circ h$, where   $t_a$ is the translation map with respect to the section
$a=f(e)$, and $h(e)=e$. We have $t_a\in G(R)$ by Corollary \ref{induced action trivial}. This reduces us to the case $f(e)=e$.
Then $f$ respects the group law, according to \cite{Mumford; Fogarty; Kirwan 1993}, Corollary 6.2.

Suppose for the moment that $R=k$ is a field.
According to \cite{Mumford 1970},   Corollary on page 132 the classifying map 
$$
P\lra \Pic^\tau_{B/k},\quad x\longmapsto [\shP|\{x\}\times B]
$$
stemming from a normalized Poincar\'e sheaf $\shP$ on $P\times B$ is an isomorphism
of abelian varieties. 
Moreover, one easily sees that it is natural with respect to $P$.
In particular, the homomorphism $f:P\ra P$ coincides with $(f^*)^*$.
Since $f^*=\id_B$, it follows that $f=\id_P$.
We refer to Section 8 for a detailed discussion  of Poincar\'e sheaves.
Moreover, the above  reasoning immediately carries  over to a noetherian base ring $R$.
%

Summing up, we have shown  that $G\subset\Aut_{P/S}$ is representable by closed embeddings,
and that the structure morphism $G\ra S$ is a family of abelian varieties.
It remains to verify that the inclusion $G\subset\Aut_{P/S}$ is open.
For this it suffices to check that for a given automorphism  $f:P\ra P$, the   set
$$
U=\{s\in S\mid \text{$f_s:P_s\ra P_s$ is a translation}\}
$$
is   open. Since this is a closed set in a noetherian topological space, the task is  to verify
that it is stable under generization.
So we may assume that $R$ is local and    $f\otimes\kappa(s)$ is a translation for the closed point $s\in S$, and have to show
that $U=S$.
Using Grothendieck's Comparison Theorem (\cite{EGA IIIa}, Theorem 5.4.1 for schemes  and \cite{SP}, 
Lemma 0A4Z for algebraic spaces), 
it suffices to treat the case that $R$ is artinian. 
Choose a composition series $\maxid_R=\ideala_0\supset\ldots\supset\ideala_r=0$ whose quotients have length one,
and write $R_i=R/\ideala_i$.
We now show by induction on $i\geq 0$ that $f_i=f\otimes R_i$ are translations. The case $i=0$ is trivial.
Suppose now $i>0$, and that $f_{i-1}$ is a translation.
Recall that $H^0(P_0,\Theta_{P_0/k})$ is the Lie algebra for the group scheme $\Aut_{P/S}\otimes k$,
where $\Theta_{P_0/k}=\uHom(\Omega^1_{P_0/k},\O_{P_0})$ is the tangent sheaf. This coincides with the Lie algebra for the group scheme of translations, and it 
follows that each  $g\in\Aut_{P/S}(k[\epsilon])$ with $g\otimes k=e\otimes k$ is a translation.

Since $P\ra S$ is smooth, we may extend $f_{i-1}\in P(R_{i-1})$ so some translation $f_i'$ over $R_i$.
Now both $f_i$ and $f_i'$  are extensions of $f_{i-1}$. 
As explained in \cite{Talpo; Vistoli 2013}, Corollary 4.4 we may view the difference as an  element in $H^0(P_0,\Theta_{P_0/k})\otimes \ideala_i/\ideala_{i-1}$,
and conclude that $f_i$   is  a translation.
\qed

\begin{proposition}
\mylabel{free and transitive}
The canonical action of $G$ on $P$ is free and transitive.
\end{proposition}

\proof
We have to show that the morphism $(\mu,\pr_2):G\times_S P\ra P\times_SP$ is an isomorphism, where
$\mu:G\times P\ra P$ denotes the action. 
It 
suffices to treat the case that $P\ra S$ is a family of abelian varieties (Proposition \ref{group law}), and that $S$ is the spectrum of 
a noetherian ring $R$.
Then we saw in the proof for Proposition \ref{automorphism group} that the $G$-action coincides with the $P$-action via
translations. The latter is obviously free and transitive.
\qed

\medskip
Let us sum up the content of the preceding propositions:

\begin{theorem}
\mylabel{para-abelian becomes torsor}
The  family of para-abelian varieties $P\ra S$ induces, in the above canonical way, a family of abelian varieties $G\ra S$
inside $\Aut_{P/S}$, and $P$ becomes a  principal homogeneous $G$-space, in other words a  representable $G$-torsor.
\end{theorem}

In turn, we get a cohomology class $[P]\in H^1(S,G)$, where the cohomology is taken with respect to
the fppf topology.  Since the structure morphism $G\ra S$ is smooth, the cohomology  remains unchanged if computed with the \'etale
topology (\cite{GB III}, Theorem 11.7). Indeed, we already saw in Corollary \ref{etale surjection} 
that $P\ra S$ admits sections locally in the \'etale topology.
Note, however,  that the order of the cohomology class $[P]$ may be infinite (\cite{Raynaud 1970}, Chapter XIII, Section 3.2
and also \cite{Zomervrucht 2015}).

Now let $f:P_1\ra P_2$ be a morphism  between families of para-abelian varieties, and $G_i\subset\Aut_{P_i/S}$
be the resulting families of abelian varieties as above, such that  $P_i$ is a principal homogeneous $G_i$-space.
 
\begin{proposition}
In the above situation, there is a unique homomorphism $f_*:G_1\ra G_2$ between families of abelian varieties
such that $f:P_1\ra P_2$ is equivariant with respect to the action of $G_1$.   
\end{proposition}

\proof
Uniqueness is clear: For each $R$-valued point $\sigma\in G_1(R)$, there is an fppf extension $R\subset R'$ such that
there is some $a'\in P_1(R')$. Then
$$
f(\sigma_{R'} \cdot a') = f_*(\sigma_{R'}) \cdot f(a') = f_*(\sigma)_{R'} \cdot f(a').
$$
Consequently $f_*(\sigma)\in G_2(R)$ is uniquely determined by $f:P_1\ra P_2$.

We now verify existence. In light of the uniqueness and fppf descent, it suffices to treat the case that there
is a section $e_1\in P_1(S)$. Composition with $f$ yields a section $e_2\in P_2(S)$.
In turn, we obtain identifications $G_i=G_i\cdot e_i= P_i$. With respect to these identifications,
we can regard $f:P_1\ra P_2$ as a morphism $f_*:G_1\ra G_2$ between families of abelian varieties
that respect  the zero sections. This is a homomorphism by \cite{Mumford; Fogarty; Kirwan 1993}, Corollary 6.4.
Given $R$-valued points $\sigma\in G_1(R)$ and $a\in P_1(R)$, we write $a=\eta\cdot e_1$ for some $\eta\in G_1(R)$
and obtain
$$
f(\sigma\cdot a) = f(\sigma\eta\cdot e_1) = f_*(\sigma\eta)\cdot e_2 = f_*(\sigma)f_*(\eta) \cdot e_2 = f_*(\sigma)\cdot f(\eta\cdot e_1) = f_*(\sigma) \cdot f(a).
$$
So $f:P_1\ra P_2$ is equivariant with respect to the $G_1$-actions stemming from  the inclusion
$G_1\subset\Aut_{P_1/S}$ and the homomorphism $ G_1\stackrel{f_*}{\ra} G_2\subset\Aut_{P_2/S}$.
\qed

\medskip
Each $R$-valued point $a\in P(R)$ can be seen as an isomorphism $\xi:G_R\ra P_R$, which comes
with an inverse $\varphi=\xi^{-1}$, and each $l\in\Pic_{G/S}(R)$ yields some $\varphi^*(l)\in \Pic_{P/S}$.
As explained in \cite{Raynaud 1970}, Chapter XIII, Proposition 1.1 this induces a  canonical isomorphism
\begin{equation}
\label{wedge identification}
\Pic_{G/S}\wedge^GP\lra \Pic_{P/S},\quad (l,a)\longmapsto \varphi^*(l).
\end{equation}
Here the wedge symbol denotes the quotient of $\Pic_{G/S}\times X$ by the diagonal left action $g\cdot(l,a)=(lg^{-1},ga)$.
The $G$-action on  $P$ is free and transitive (Proposition \ref{free and transitive}),
whereas the $G$-action on the  invariant open subspace $\Pic^\tau_{G/S}$ is   trivial (Corollary \ref{induced action trivial}).
It follows that the projection 
\begin{equation}
\label{projection bijective}
\pr_1:\Pic^\tau_{G/S}\wedge^GP\lra (\Pic^\tau_{G/S})/G=\Pic^\tau_{G/S}
\end{equation} 
is an isomorphism. As observed  by Raynaud in loc.\ cit., composing \eqref{wedge identification} with the inverse of \eqref{projection bijective} 
yields:

\begin{proposition}
The above maps gives  an identification $\Pic^\tau_{G/S}=\Pic^\tau_{P/S}$ of families of abelian varieties.
\end{proposition}

For each family $A\ra S$ of abelian varieties, the   family of abelian varieties $\Pic^\tau_{A/S}$
is called the \emph{dual family}. 
The above observation identifies  the family $B=\Pic^\tau_{P/S}$ of abelian varieties  coming
from our family $P\ra S$ of para-abelian varieties  with the dual family for $G\ra S$, where  the latter is
defined via the inclusion $G\subset\Aut_{P/S}$.

This has a remarkable consequence:
Let $P\ra S$ and $P'\ra S$ be two families of  para-abelian varieties, and consider the canonical map
$$
\Hom_S(P',P)\lra \Hom_{\text{Gr}/S}(\Pic^\tau_{P/S},\Pic^\tau_{P'/S}),\quad f\longmapsto f^*.
$$
Note that the term on the left is   a set, which might be empty,
whereas the term on the right is an abelian group. Moreover, the $G$-action   on $P$ induces an action of the group  $G(S)$ on the set $\Hom_S(P',P)$.

\begin{lemma}
\mylabel{orbits in hom-set}
In the above situation, the action of    $G(S)$ on the set $\Hom_S(P',P)$ is free,
and the   fibers of the map $f\mapsto f^*$ are precisely the orbits.
\end{lemma}

\proof
Suppose some $\sigma\in G(S)$ fixes a morphism $f:P'\ra P$.
We have to verify that $\sigma=e$. By descent, we may replace
$S$ with $P'$ and assume that there is a section $a\in P'(S)$.
Then $f(a)\in P(S)$.  The action of $G(S)$ on $P(S)$ is free, and  from $\sigma + f(a) = f(a)$ it follows $\sigma=e$.

Corollary \ref{induced action trivial} ensures that the orbits are contained in the fibers.
Conversely, suppose that $f,g:P'\ra P$ are two morphisms with $f^*=g^*$.
We have to produce some $\sigma\in G(S)$ with $f=g+\sigma$.
It must be unique, if it exists, according to the preceding paragraph.
By descent, our problem is local, so we     may assume that $P'\ra S$ and hence also $P\ra S$ is a family of abelian varieties,
that $S$ is the spectrum of a noetherian ring $R$, and that we  have  an identification $G=P$.
Write $f= t_a\circ f_0$ and $g=t_b\circ g_0$, where $f_0,g_0:P'\ra P$ are homomorphisms,
and $t_a,t_b$ are translations by some elements $a,b\in P(S)=G(S)$.
Our task is to verify that $g_0=f_0$. In light of 
Grothendieck's Comparison Theorem (\cite{EGA IIIa}, Theorem 5.4.1 for schemes  and \cite{SP},  
Lemma 0A4Z for algebraic spaces)
it suffices to check
this if $R$ is a local Artin ring.
From Corollary \ref{induced action trivial} we have  $f_0^*=g_0^*$.
We now argue as in  Proposition \ref{automorphism group}, by choosing a composition series
$\maxid_R=\ideala_0\supset \ldots\supset\ideala_r=0$ whose quotients have length one,
and applying  induction on $i\geq 0$ with $R_i=R/\ideala_i$.
\qed

\section{Some extensions of group schemes}
\mylabel{Some extensions}

Let $k$ be a ground field of characteristic $p\geq 0$. 
Recall that an \emph{abelian variety} $A$ is a group scheme that is smooth, connected and proper.
By abuse of notation, we simply say that $A$ is \emph{abelian}.
A group scheme   $N$ is called \emph{of multiplicative type} if there is a field extension $k\subset k'$
such that $N\otimes k'$ is isomorphic to the spectrum of the Hopf algebra $k'[\Lambda]$ for some
commutative group $\Lambda$. We are mainly interested in the case that the scheme $N$ is of finite type; then
the group $\Lambda$ is finitely generated, and one may choose $k\subset k'$ finite and separable.
Moreover, $N$ is a twisted form, already in the \'etale topology, of  
$\GG_m[n_1]\oplus\ldots\oplus\GG_m[n_r]$, with certain invariant factors $n_r|\ldots|n_1$.
Here the summands are the kernels of the multiplicative group $\GG_m$
with respect to multiplication by $n_i\geq 0$. For brevity, we call such $N$   \emph{multiplicative}.

Throughout this section, we study   extensions of group schemes
\begin{equation}
\label{extension}
0\lra N\lra E\lra A\lra 0,
\end{equation}
where $A$ is abelian,  $N$ is multiplicative, and the middle term $E$ is commutative, and analyze their splittings.
Some of the assertions below are valid over general base schemes, but for the sake of exposition we stick to a ground field $k$.
First note that $\Hom(A,N)=0$, because $N$ is affine and $h^0(\O_A)=1$.
Hence a splitting is unique, if it exists.

Recall that \emph{exactness} means that for each $T=\Spec(R)$,
the sequence of groups $0\ra N(R)\ra E(R)\ra A(R)$ is exact, and for each  $a\in A(R)$ there is an fppf extension $R\subset R'$
such that the base-change $a'$ is in the image of $E(R')\ra A(R')$.
 
For any commutative group scheme   $G$, we can consider the resulting 
abelian sheaves $\uExt^i(G,\GG_m)$ on the category $(\Aff/k)$ of affine schemes,
endowed with the fppf topology. These are defined with injective resolutions of $\GG_m$,
but can also be seen as sheafifications of the presheaves that assign to each ring $R$
the groups $\Ext^i(G_R,\GG_{m,R})$. Its elements  can also be interpreted as equivalence classes of Yoneda extensions,
formed with sheaves of abelian groups.
Note that the $\uExt^i(G,\GG_m)$  may or may not be representable by group schemes. Let us recall the following facts:

\begin{proposition} 
\mylabel{zero sheaves}
We have $\uExt^1(N,\GG_m)=0$.  
Moreover,  $\uHom(A,H)=0$ for any affine group scheme $H$. 
\end{proposition}

\proof 
The first is contained in \cite{SGA 7a}, Expos\'e VIII, Proposition 3.3.1. 
For the second assertion, suppose  $f:A_R\ra H_R$ is a homomorphism over some ring $R$.
Since the formation of $\Gamma(A,\O_A)=k$ commutes with flat ring extensions, we have $R=\Gamma(A_R,\O_{A_R})$,
and conclude that $f$ factors over the zero section.
\qed
 
\medskip 
We infer that the short exact sequence \eqref{extension} yields an exact sequence
$$
0\lra \uHom(E,\GG_m)\lra \uHom(N,\GG_m)\lra \uExt^1(A,\GG_m)\lra \uExt^1(E,\GG_m)\lra 0
$$
of abelian sheaves. The sheaf $\uExt^1(A,\GG_m)$ is representable by an abelian variety,
and the theory of bi-extension gives an identification with the \emph{dual abelian variety} $B=\Pic^\tau_{A/k}$,
see the discussion in \cite{Roessler; Schroeer 2020}, Section  2. Actually, $A\mapsto B$ is an anti-equivalence of
the category of abelian varieties with itself, coming with the \emph{biduality identification} $A=\uExt^1(B,\GG_m)$.

Let us call a group scheme $L$ a \emph{local system} if it is a twisted form of the constant group scheme
$(\Lambda)_k$, where $\Lambda$ is a finitely generated abelian group.
Equivalently, $L$ is \'etale, and the group $L(k^\sep)$ is finitely generated.
For each multiplicative group scheme $G$, the sheaf $L=\uHom(G,\GG_m)$ is a local system.
Actually, the functor $G\mapsto L$ is an anti-equivalence between the   category of multiplicative group schemes  
and the category of local systems, which follows from \cite{SGA 3b}, Expos\'e IX, Corollary 1.2.
Summing up, our extension \eqref{extension} gives a coboundary map
$$
L=\uHom(N,\GG_m)\stackrel{\partial}{\lra} \uExt^1(A,\GG_m)=B
$$
from the local system $L$ to the dual abelian variety $B$.
 
\begin{proposition}
\mylabel{coboundary}
Our  extension \eqref{extension} splits if and only if   $\partial:L\ra B$ vanishes.
\end{proposition}

\proof
The condition is obviously necessary. Suppose now that the coboundary map vanishes. 
We already remarked that the section is unique, if it exists,  so with Galois descent it suffices
to treat the case that $k$ is separably closed. Then $N$ is a direct sum of copies of the multiplicative group $\GG_m$ and the kernels $\mu_n=\GG_m[n]$,
hence it is enough to consider the cases $N=\GG_m$ and $N=\mu_n$.
In the former case, the extension class of \eqref{extension} is the image of the identity map $\id:N\ra \GG_m$ under the coboundary map,
whence the extension splits.
It remains to treat the case $N=\mu_n$. The canonical inclusion $\mu_n\subset\GG_m$ yields
a push-out extension $E'=(E\oplus\GG_m)/\mu_n$, which yields a commutative diagram
$$
\begin{CD}
\uHom(\GG_m,\GG_m) @>\partial>> \uExt^1(A,\GG_m)\\
@VVV                    @VV\id V\\
\uHom(\mu_n,\GG_m) @>>\partial> \uExt^1(A,\GG_m).
\end{CD}
$$
On the right we have the dual abelian variety $B=\uExt^1(A,\GG_m)$.
The image of the identity on $\GG_m$ is a rational point $b\in B$, and we have $b=0$ by assumption.
  Thus $E'$ splits. The Kummer sequence gives an exact sequence
$$
\Hom(A,\GG_m)\lra \Ext^1(A,\mu_n)\lra \Ext^1(A,\GG_m).
$$ 
The term on the left vanishes. Thus the map on the right is injective, and $E$ splits as well.
\qed

\medskip
The coboundary map is functorial in the extension \eqref{extension}, by the very definition of delta functors
(\cite{Grothendieck 1957}, Section 2.1).
Consequently, the pull-back along some homomorphism of abelian varieties $A'\ra A$ splits if and only
if the composition $L\ra B\ra B'$ vanishes.
We exploit this as follows:
Let $Z\subset B$ be the Zariski closure of the set-theoretical image for $L\ra B$.
 This is a smooth subgroup scheme, and its formation commutes with ground field extensions
(\cite{SGA 3b}, Expos\'e VIb, Proposition 7.1).
The short exact sequence $0\ra Z\ra B\ra B'\ra 0$ defines an abelian variety $B'=B/Z$, and we now consider
the dual abelian variety $A'$. It comes with a homomorphism
$$
A'=\uExt^1(B',\GG_m)\stackrel{f}{\lra} \uExt^1(B,\GG_m)=A.
$$

\begin{theorem}
\mylabel{splitting map}
The homomorphism $f:A'\ra A$ has the following properties:
\begin{enumerate}
\item The formation of $A'$ and $f$ commutes with ground field extensions.
\item The pull-back of our extension \eqref{extension} along $f:A'\ra A$ splits.
\item The induced homomorphism $A'\ra E$ is a closed embedding, and its image contains every abelian subvariety inside $E$.
\item For every homomorphism $g:A''\ra A$ from an abelian variety $A''$ such that the pullback of \eqref{extension} along $g$ splits,
there is a unique factorization over $f:A'\ra A$.
\end{enumerate}
\end{theorem}
 
\proof
Assertion (i) follows from the fact that the formation of $B'$ commutes with ground field extensions, 
whereas (ii) is a consequence of Proposition \ref{coboundary}.
We next verify (iv). Suppose our extension is split by some $g:A''\ra A$. Let $B\ra B''$ be the dual homomorphism.
Its kernel $K\subset B$ is a closed subscheme. We already observed that the composite map $L\ra B\ra B''$ vanishes.
Thus the set-theoretical image of $L\ra B$ is contained in $K$, hence the Zariski closure $Z$ is contained
in $K$. The Isomorphism Theorem gives a unique factorization of $B\ra B'\ra B''$.
Dualizing gives the desired factorization $A''\ra A'\ra A$. This factorization is unique, by biduality.

It remains to establish (iii). The splitting for $E'=E\times_AA'$ is unique, as we already observed below \eqref{extension}.
In turn, there is a unique lift $g:A'\ra E$ for $f:A'\ra A$. Consider the schematic image $A_0=g(A')$ inside $E$,
and the induced factorization $A'\ra A_0\ra A$, and the   dual factorization $B\ra B_0\ra B'$.
Obviously, the  pull-back of \eqref{extension} along $A_0\ra A$ splits. From (iv) we get a decomposition
$A'=A_0\oplus A_1$. This ensures that  $B_0\ra B'$ is a closed embedding. Using that $B\ra B'$ is surjective,
we infer that $B_0=B'$ is bijective, and it follows that $A'\ra A_0$ is an isomorphism.
\qed

\medskip
The sheaf kernel $\Kernel(\partial)$ and the sheaf image $\Image(\partial)=L/\Kernel(\partial)$ for the homomorphism   
$\partial:L\ra B$ are group schemes that are locally of finite type, and in fact local systems. 
According to Lemma \ref{quotient}, the sheaf cokernel $\Cokernel(\partial)=B/\Image(\partial)$
is an algebraic space that is not necessarily locally separated.

\section{Maximal abelian subvarieties}
\mylabel{Maximal abelian}

Let $k$ be a ground field of characteristic $p\geq 0$, and $G$ be a group scheme of finite type.
The goal of this section is to describe the  maximal abelian subvariety $G_\ab\subset G$, making
evident that its formation commutes with ground field extensions. 
We start by defining a   three-step filtration $G=G_0\supset G_1\supset G_2\supset G_3$, which is of independent interest.
 
Let $G^\aff$ be the spectrum of the ring $\Gamma(G,\O_G)$.
Then the group law on $G$ induces a group law on $G^\aff$, and the canonical map $G\ra G^\aff$ is a homomorphism.
We denote by $G_1$ its kernel. Then $G\ra G^\aff$ is flat and surjective, so that $G/G_1=G^\aff$, and furthermore $h^0(\O_{G_1})=1$
(\cite{Demazure; Gabriel 1970}, Chapter III, \S3, Theorem 8.2).  
The latter condition means that $G_1$ is \emph{anti-affine}. 
According to \cite{Brion 2017}, Proposition 3.3.4 this ensures that $G_1$ is smooth and commutative.
We then define $G_2\subset G_1$ as the largest subgroup scheme that is smooth, connected and affine (\cite{Brion 2017}, Lemma 3.1.4).
Finally, write $G_3\subset G_2$ for the largest subgroup scheme that is multiplicative, so that $G_2/G_3$
is unipotent (\cite{Demazure; Gabriel 1970}, Chapter IV, \S3, Theorem 1.1). This defines the desired  three-step   filtration on $G$.
 
\begin{proposition}
\mylabel{three step filtration}
The three-step filtration on $G$ has the following properties:
\begin{enumerate}
\item Each homomorphism $f:G\ra G'$ respects the filtrations.
\item The formation of $G_i\subset G$ commutes with ground field extensions.
\item The group scheme $G_1$ is anti-affine, and    $G_1/G_2$ is an abelian variety.
\item The extension $0\ra G_3\ra G_2\ra G_2/G_3\ra 0$ has a unique splitting.
 \item We have $G_2=G_3$   in characteristic $p>0$, whereas $G_2/G_3\simeq\GG_a^{\oplus r}$ for some $r\geq 0$
in characteristic zero.
\end{enumerate}
\end{proposition}

\proof
We start with assertion (iii). We already observed above that $G_1$ is anti-affine.
In characteristic $p>0$ it is also semi-abelian, so that there is a short exact sequence
\begin{equation}
\label{semiabelian}
0\lra T\lra G_1\lra A\lra 0,
\end{equation}
for some torus $T$ and some abelian variety $A$, according to \cite{Brion 2009}, Proposition 2.2.
The torus is contained in $G_2$, by maximality of the latter. The resulting surjection $G_1/T\ra G_1/G_2$ reveals that $G_1/G_2$ is abelian.
In characteristic zero we use  the smallest subgroup scheme $N\subset G_1$ such that the quotient  $G_1/N$ is proper.
This exists in all characteristics, and is affine and connected, according to \cite{Brion 2017}, Theorem 2.
For $p=0$ the group schemes $G_1/N$ and $N$ are automatically smooth. Consequently 
$G_1/N$ is abelian and $N\subset G_2$, and  we conclude again that $G_1/G_2$ is abelian.

Next we consider assertion (i). Let $f:G\ra G'$ be a homomorphism. The composite morphism $G_1\ra G'^\aff$ vanishes, because $h^0(\O_{G_1})=1$.
The next composite   $G_2\ra G'_1/G'_2$ also vanishes: Its image is a quotient of $G_2$ and thus
smooth, connected and affine. It is also  a closed subgroup scheme in $G'_1/G'_2$, hence proper, and therefor trivial.
Finally, $G_3\ra G'_2/G'_3$ is zero, because the domain is multiplicative and the range is unipotent 
(\cite{Demazure; Gabriel 1970}, Chapter IV, \S 3, Proposition 1.3).

We now come to (iv). The extension in question has at most one splitting, because  there are no non-trivial homomorphisms from the unipotent
group scheme $G_2/G_3$ to the multiplicative group scheme $G_3$.
If $k$ is perfect, such a splitting indeed exists, by \cite{Demazure; Gabriel 1970}, Chapter IV, \S 3, Theorem 1.1.
This also ensures that $G_3$ is smooth and connected.
Suppose now $p>0$. Then we have an exact sequence \eqref{semiabelian}.
The canonical projection $G_3\ra A$ vanishes, and we obtain a commutative diagram
$$
\begin{CD}
0@>>> G_3 @>>> G_2 @>>> G_2/G_3 @>>> 0\\
@.      @VVV    @VVV    @VVV\\
0@>>> T @>>>    G_1 @>>> A@>>>0.
\end{CD}
$$
The induced map $G_2/G_3\ra A$ is zero, so the Snake Lemma gives an inclusion $G_2/G_3\subset T/G_3$.
Consequently,  $G_2/G_3$ is both unipotent and  multiplicative,
hence zero. In turn, the extension splits for trivial reasons.
This also establishes the first part of (iv). In characteristic zero, the   unipotent scheme $G_2/G_3$ must be   smooth. 
It follows from \cite{Demazure; Gabriel 1970}, Chapter IV, \S 2 that it is isomorphic to a sum of additive groups $\GG_a$.

It remains to establish (ii).
For each quasicompact and quasiseparated scheme $X$, the formation of $H^0(X,\O_X)$ commutes with ground field extensions $k\subset k'$.
In turn, the same holds for the kernel $G_1\subset G$.
Suppose $N'$ is a smooth connected affine subgroup scheme in the   base-change $G'_1=G_1\otimes k'$ that contains 
$G_2'=G_2\otimes k'$. Then the quotient $H'=N'/G_2'$ is smooth, connected and affine,
and comes with an embedding into the base-change of the abelian variety $A=G_1/G_2$. It follows that $H'=0$.
Thus $G_2\subset G_1$ commutes with base-change.
Finally, the inclusion $G_3\subset G_2 $ commutes with base-change by \cite{Demazure; Gabriel 1970}, Chapter IV, \S 3,  Proposition 1.3.
\qed

\medskip
Let $H\subset G$ be an abelian subvariety of largest dimension.
Then $H\subset G_1$, because  $G/G_1$ is affine and $h^0(\O_H)=1$.
Since $G_1$ is commutative, and    sums and quotients of abelian varieties remain abelian,
every other abelian subvariety $H'\subset G_1$ must be contained in $H$.
Therefore, $G_\ab=H$ is the largest abelian subvariety, and  we also call it 
the \emph{maximal abelian subvariety} of $G$. Clearly, this is functorial in $G$.
However, it is not immediately evident that its formation commutes with ground field extensions.

In characteristic $p>0$, we have $G_2=G_3$ and   an extension $0\ra T\ra G_1\ra A\ra 0$
of the abelian variety $A=G_1/G_2$ by the torus $T=G_2$.  Let $A'\subset G_1$ be the abelian subvariety
constructed in Section \ref{Some extensions} as the dual of $B'=B/Z$, and regard $A'$
as a subgroup scheme of $G$.  

\begin{theorem}
\mylabel{maximal abelian base-change}
The formation of the maximal abelian subvariety $G_\ab\subset G$ commutes with ground field extensions $k\subset k'$.
In characteristic $p>0$  we have $G_\ab=A'$.
\end{theorem}

\proof
Let $k\subset\Omega$ be a field extension and write $G_\Omega=G\otimes \Omega$.
It suffices to check that the maximal abelian subvarieties inside $G$ and $G_\Omega$ have the same dimension.
Seeking a contradiction, we assume that there is an abelian subvariety $N\subset G_\Omega$ with $\dim(N)>\dim(G_\ab)$.
Replacing $G$ by the subgroup scheme $G_1$ we reduce to the case that $G$ is commutative,
hence  $G_\ab$ must be normal. By passing to $G/G_\ab$ we may also assume that $G_\ab=0$.
We now have a non-trivial abelian subvariety $N\subset G_\Omega$, and will reach a contradiction by producing
a non-trivial abelian subvariety $H\subset G$.  Note that for along the way  we may enlarge $\Omega$.

We now argue as in \cite{Schroeer 2022}, proof for Theorem 6.1:
Using \cite{EGA IVc}, Theorem 8.8.2 we   reduce to the case that  the field extension $k\subset\Omega$ is finitely generated.
By considering suitable intermediate field and enlarging  $\Omega$ if necessary, it suffices to treat the 
cases that $k\subset\Omega$ is either purely transcendental, or a finite Galois extension, or a finite radical extension in characteristic $p>0$.
In the first case, we extend $N\subset G_\Omega$ to a family of abelian varieties in $G_R$ over some
localization $R=k[T_1,\ldots,T_n]_f$, and obtain a contradiction by specializing to a rational point
in $\Spec(R)$. 
In the second case we use that the maximal abelian subvariety of $G_\Omega$ is stabilized by the elements
of the Galois group $\Gamma=\Gal(\Omega/k)$, hence descends to an abelian subvariety in $G$.
In the third case we are in characteristic $p>0$, and this was 
 essentially solved above:  From Theorem \ref{splitting map} we   get  $G_\ab=A'$, and its formation   commutes with ground field
extensions. 
\qed

\section{The notion of Albanese maps}
\mylabel{Albanese maps}
 
Let $S$ be a base scheme, and $X$ be an algebraic space where the structure morphism $X\ra S$
is proper, flat, of finite presentation, and  cohomologically flat in degree $d=0$. 
Then $\Pic_{X/S}^\tau$ exists as an algebraic space whose structure morphism is of finite presentation (hence quasiseparated) and locally separated, according to Theorem \ref{pic tau}.
We now come to the central topic of this paper, where we re-define and generalize 
classical notions of Albanese varieties and Albanese maps:

\begin{definition}
\mylabel{definition albanese map}
An \emph{Albanese map} for $X$ is a pair $(P,f)$ where $P\ra S$ is a family of para-abelian varieties 
and $f:X\ra P$ is a morphism satisfying the  following condition: For each   $s\in S$ the homomorphism $f^*$ identifies
$\Pic^\tau_{P/S}\otimes \kappa(s)$ with the maximal abelian subvariety of $\Pic^\tau_{X/S}\otimes\kappa(s)$.
\end{definition}

By abuse of notation, \emph{we simply say that $f:X\ra P$ is an Albanese map}.
From Lemma \ref{kernel trivial} we see that  $f^*:\Pic^\tau_{P/S}\ra \Pic^\tau_{X/S}$ is a monomorphism.
It is actually a  closed embedding provided that $\Pic^\tau_{X/S}$ is separated, and the latter indeed
holds if $S$ is artinian.  
 
A priori, our notion of Albanese maps has good base-change properties:
Let $f:X\ra P$ be   any morphism to a family $P$ of para-abelian varieties. Given  $S'\ra S$ we write
$X'=X\times_SS'$ and $P'=P\times_SS'$ for the base-changes, and $f':X'\ra S'$ for the induced morphism.

\begin{proposition}
\mylabel{base-change albanese map}
In the above situation, if $f:X\ra P$ is an Albanese map, the same holds for the base-change $f':X'\ra P'$.
The converse remains true  if   $S'\ra S$ is surjective.
\end{proposition}

\proof
This follows from the fact that the  maximal abelian subvarieties
are stable under field extensions (Theorem \ref{maximal abelian base-change}).
\qed

\medskip
The next observation is an important intermediate step towards existence, uniqueness and universal property of Albanese maps:

\begin{proposition}
\mylabel{weak uniqueness}
Suppose $f:X\ra P$ is an Albanese map. Let $g:X\ra Q$ be a morphism to some other family $Q$ of para-abelian varieties.
Then the homomorphism $g^*:\Pic^\tau_{Q/S}\ra\Pic^\tau_{X/S}$ admits  a unique factorization over $f^*:\Pic^\tau_{P/S}\ra\Pic^\tau_{X/S}$.
Moreover, there is at most one morphism  $h:P\ra Q$ such that the diagram
\begin{equation}
\label{maps to para-abelian}
\begin{tikzcd} 
                & X\ar[dl,"f"']\ar[dr,"g"]\\
P\ar[rr,"h"']    &                           &Q
\end{tikzcd}
\end{equation}
is commutative.
\end{proposition}

\proof 
Since $f^*$ is a monomorphism, there is at most one  factorization for $g^*$.  
Consequently, the existence of a factorization is   local in $S$ for the fpqc topology,  and it suffices to treat the
case that $S$ is the spectrum of a ring $R$. Since $X$ is of finite presentation, we moreover
may assume that $R$ is noetherian.
Set
$$
G=\Pic^\tau_{X/S}\quadand N=\Pic^\tau_{P/S}\quadand H=\Pic^\tau_{Q/S}.
$$
The structure morphism  $N\ra S$ is fppf, and the  translation action of $N$
on $G$ is free, hence the quotient $G/N$ exists as an algebraic space, and its formation commutes with base-change,
according to Lemma  \ref{quotient}.
We see that $g^*:H\ra G$ factors over $N$ if and only
if the composite homomorphism $c:H\ra G/N$ is trivial.  
Note that the structure morphism $G/N\ra S$ is of finite type, because the same
holds for $G\ra S$  and the projection $G\ra G/N$ is   fppf.
Moreover, for each point $s\in S$, the inclusion $N_s\subset G_s$ is the maximal abelian subvariety, so the quotient $(G/N)_s=G_s/N_s$ does not contain
any  non-zero abelian subvarieties.

Now suppose that $R$ is a    local Artin ring, such that $S=\{s\}$. The closed fiber  $(G/N)_s$ is a separated scheme. Hence also $G/N$ is separated,
and it is schematic by \cite{Rydh 2015},  Corollary 8.2. 
The set-theoretic image of  $H\ra G/N$ is the origin   $e_s\in (G/N)_s$ in the closed fiber.
By the Rigidity Lemma (\cite{Mumford; Fogarty; Kirwan 1993}, Proposition 6.1), the morphism $H\ra G/N$ must factor over some section $\sigma:S\ra G/N$.
We must have $\sigma=e$, because the image  of the map  $H(R)\ra G/N(R)$ is a  subgroup. 

Using Grothendieck's Comparison Theorem (\cite{EGA IIIa}, Theorem 5.4.1 for schemes and \cite{SP},
Lemma 0A4Z for algebraic spaces), 
we infer that $H\ra G/N$ is trivial if $R$ is a     local
noetherian ring that is complete. For general local noetherian rings  $R$, the formal completion $ \widehat{R}=\invlim R/\maxid_R^n$ is an fpqc extension,
and it follows from uniqueness applied over $\widehat{R}\otimes_R\widehat{R}$ and  fpqc descent that $H\ra G/N$ then is trivial as well.
For general noetherian rings $R$, one sees that for each prime $\primid\subset R$, there is an element $f\in R\smallsetminus\primid$
such that $H\ra G/N$ is trivial over the localization $R_f$. Hence $H\ra G/N$ is trivial over the whole spectrum $S=\Spec(R)$.

Suppose there are two morphisms $h_1,h_2:P\ra Q$ making the above diagram commutative.
It remains to show that $h_1=h_2$. We just saw  that the induced homomorphisms
$h_i^*:\Pic^\tau_{Q/S}\ra\Pic^\tau_{P/S}$ coincide. Let $I\subset\Aut_{Q/S}$ be the inertia subgroup scheme
for the numerically trivial sheaves. We showed  in Section \ref{Equivariance} that this is a family of abelian varieties,
and its action on $Q$ is free and transitive. By Lemma \ref{orbits in hom-set} there is a unique section $\sigma\in I(S)$
with $h_2=\sigma+h_1$.
By fppf descent, we may replace $S$ with $X$ and assume that the structure morphism $X\ra S$
admits a section $\tau\in X(S)$. The commutativity of the diagram \eqref{maps to para-abelian} reveals that the 
$h_1$ and $h_2$ coincide on the section $f\circ \tau\in P(S)$, thus $\sigma\in I(S)$ must be trivial.
This shows $h_1=h_2$.
\qed

\medskip
In the situation of the Proposition, let  $S'\ra S$ be some   morphism of schemes, and write $X',P',Q'$
for the base-changes, and $f':X'\ra P'$  and $g':X'\ra Q'$ for the induced maps.

\begin{corollary}
\mylabel{morphism descends}
Suppose that $S'\ra S$ is fpqc, and that there is some $h':P'\ra Q'$  with $g'=h'\circ f'$. Then there is also a morphism $h:P\ra Q$ with $g=h\circ f$,
and $h'$ equals the base-change of $h$.
\end{corollary}

\proof
Consider the fiber product $S''=S'\times_SS'$, which comes with two projections $\pr_i:S''\ra S'$. By the uniqueness in the proposition applied over $S''$,
we have an equality $\pr_1^*(h')=\pr_2^*(h')$. From fpqc descent 
(\cite{SGA 1}, Expos\'e VIII, Theorem 5.2 for schemes and \cite{SP}, 
Lemma 0ADV for algebraic spaces)
one deduces that there is a unique $S$-morphism $h:P\ra Q$ inducing the $S'$-morphism $h':P'\ra Q'$.
The property $g=h\circ f$ can be checked after base-changing to $S'$, where it holds by assumption.
\qed

\section{Poincar\'e sheaves}
\mylabel{Poincare}

We now review the notion of Poincar\'e sheaves, which will be used in the next section to construct Albanese maps.
Fix a base scheme $S$, and let $X$ be an algebraic space whose structure morphism $\varphi:X\ra S$
is proper, flat, of finite presentation, and cohomologically flat in degree $d=0$, with $h^0(\O_{X_s})=1$ 
for all points $s\in S$.
 Thus $\varphi_*(\O_X)=\O_S$, and this commutes with base-change.
Consequently   $\varphi_*(\GG_{m,X})=\GG_{m,S}$, and this also commutes with base-change.
By the very definition, the algebraic space  $\Pic_{X/S}$ represents the higher direct image   $R^1\varphi_*(\GG_{m,X})$.
For each $T$, the  Leray--Serre spectral sequence for the projection $\pr_2:X_T\ra T$, together with Hilbert 90,  gives an exact five-term sequence
\begin{equation}
\label{leray serre}
0\lra \Pic(T)\lra \Pic(X_T)\lra \Pic_{X/S}(T)\lra H^2(T,\GG_m) \stackrel{\pr_2^*}{\lra} H^2(X_T,\GG_m).
\end{equation}
We are particularly interested in the case  that $T$ coincides with $\Pic_{X/S}$, or at least comes
with a morphism to $\Pic_{X/S}$:

\begin{definition}
\mylabel{poincare sheaf}
Suppose we have a morphism $f:T\ra\Pic_{X/S}$.
An invertible sheaf $\shP$ on $X\times T$ is called a \emph{Poincar\'e sheaf  with respect to $f:T\ra\Pic_{X/S}$}
if its   class in $\Pic(X_T)$ maps to the element $f\in \Pic_{X/S}(T)$  in the above sequence.
\end{definition}

By abuse of notation, we  simply   say that $\shP$ is a Poincar\'e sheaf  on $X\times T$.
Note that as customary, the product is formed over the base scheme $S$, such that $X\times T=X\times_ST$.
Poincar\'e sheaves exist  if the map $H^2(T,\GG_m) \ra H^2(X_T,\GG_m)$ 
is injective, by exactness of \eqref{leray serre}. This obviously holds  if the structure map $\varphi:X\ra S$ has a section,
or if the cohomology group $H^2(T,\GG_m)$ vanishes.
In any case, Poincar\'e sheaves are unique up to preimages of  invertible sheaves   on $T$.

Now let $\shL$ be any invertible sheaf on the product $X\times T$, for some algebraic space $T$.
This gives an element in $\Pic(X_T)$, which induces a $T$-valued point for   $\Pic_{X/S}$.
Let $f:T\ra\Pic_{X/S}$ be the resulting classifying map. Set $P=\Pic_{X/S}$, and assume there is a   Poincar\'e sheaf $\shP$  on $X\times P$
with respect to the identity map $P\ra\Pic_{X/S}$.

\begin{proposition}
\mylabel{unique pullback}
Assumptions as above.
Up to isomorphism,   there is a unique invertible sheaf $\shN$ on $T$ such that
$(\id_X\times f)^*(\shP)$ is isomorphic to $\shL\otimes \pr_2^*(\shN)$.
\end{proposition}

\proof
The five-term sequence \eqref{leray serre} is natural in $T$, so we get a commutative diagram
$$
\begin{CD}
0	@>>>	\Pic(P)	@>>>	\Pic(X_P)	@>>>	\Pic_{X/S}(P)\\
@.		@Vf^*VV		@V(\id_X\times f)^*V V	@VV\Pic_{X/S}(f)V\\
0	@>>>	\Pic(T)	@>>>	\Pic(X_T)	@>>>	\Pic_{X/S}(T).
\end{CD}
$$
By construction, both $\shL$ and $(\id_X\times f)^*(\shP)$ are Poincar\'e sheaves with respect to $f:T\ra\Pic_{X/S}$.
In light of the lower exact sequence, there class  differs by a unique element form $\Pic(T)$.
\qed

\medskip
In particular, for every field $k$ and every $l:\Spec(k)\ra T$, the isomorphism class of the pull-back $\shP|X\otimes k$
corresponds to the induced $k$-valued point $l\in \Pic_{X/S}(k)$.
We also see that in general, Poincar\'e sheaves do not exist: For example, a smooth curve $C$ of genus $g=0$
over a ground field $k$
has constant Picard scheme $\Pic_{C/k}=(\ZZ)_k$. By Riemann--Roch,
if there is a  Poincar\'e sheaf on $T=\Pic_{C/k}$, or even on $T=\{1\}$
then $C$ contains a rational point, and is thus isomorphic to the projective line.
This leads to the following criterion:

\begin{proposition}
\mylabel{for poincare}
Suppose there are morphisms $g_i:Z_i\ra X $, $1\leq i\leq r$ such that each structure morphism $h_i:Z_i\ra S$
is locally free of degree $d_i\geq 1$, with $\gcd(d_1,\ldots,d_r)=1$.
Then there is a Poincar\'e sheaf   on $X\times  \Pic_{X/S}$.
\end{proposition}

\proof
Set $T=\Pic_{X/S}$. Fix an index $1\leq i\leq r$, write  $Z=Z_i$ and $g=g_i$.
The map   on the right in \eqref{leray serre} sits in a sequence
$$
H^2(T,\GG_m)\stackrel{\pr_2^*}{\lra} H^2(X_T,\GG_m)\stackrel{g^*}{\lra} H^2(Z_T,\GG_m).
$$
Write $\psi=\pr_2\circ g$ for the composite morphism $Z_T\ra T$.
One easily  checks that the direct image sheaves $R^i\psi_*(\GG_{m,Z_T})$ vanish for all degrees $i\geq 1$ 
(compare \cite{Schroeer 2020}, proof for Lemma 1.4). 
Consequently, the Leray--Serre spectral sequence   gives an identification 
$H^2(Z_T,\GG_m)=H^2(T,\psi_*(\GG_{m,Z_T}))$. 
The composition of the inclusion $\GG_{m,T}\subset \psi_*(\GG_{m,Z_T})$ with the norm map
$N:\psi_*(\GG_{m,Z_T})\ra\GG_{m,T}$ is multiplication by $d=d_i$.
So for each $\alpha\in H^2(T,\GG_m)$ we have 
$d\cdot\alpha= N_*(g^*(\pr_2^*(\alpha)))$.
We conclude that the kernel for the map on the right in \eqref{leray serre} is annihilated by $\gcd(d_1,\ldots,d_r)=1$.
Thus $\Pic(X_T)\ra\Pic_{X/S}(T)$ is surjective, hence a Poincar\'e sheaf exists. 
\qed

\medskip 
We also have a non-existence result: 

\begin{proposition}
\mylabel{against poincare}
Suppose $S$ is the spectrum of a field $k$, and that there is a quasicompact scheme $Z$ 
so that the     pull-back map   $H^2(Z,\GG_m)\ra H^2(X_Z,\GG_m)$ is not injective.
Then there is some point $l\in\Num_{X/k}$ and some connected component
$T\subset \Pic^l_{X/k}$ such that there is no Poincar\'e sheaf with respect to $T$. 
\end{proposition}

\proof
Choose some section $s_Z\in H^0(Z,R^1\pr_{2,*}(\GG_{m,X_Z}))$ mapping to some non-trivial
element of the kernel for $H^2(Z,\GG_m)\ra H^2(X_Z,\GG_m)$. The resulting classifying map  $h:Z\ra \Pic_{X/k}$
factors over a finite union of connected components, corresponding to points $l_1,\ldots,l_r\in \Num_{X/k}$.
Moreover, $s_Z$ is the pull-back of the universal section.
The decomposition of $\Pic_{X/k}$ into connected components
gives a decomposition of the quasicompact scheme $Z$ into open-and-closed subschemes $Z_1,\ldots,Z_r$. By passing to  
one of them, we may assume $r=1$ and  let $T\subset \Pic^{l_1}_{X/k}$ be the connected component containing the image. 
Suppose that there is a Poincar\'e sheaf 
$\shP_{X\times T}$. Then the restriction $s_T$ of the universal section maps to zero under
the coboundary. In turn, its pull-back $s_Z$ also maps to zero under the coboundary, contradiction.
\qed

\section{Existence and universal property}
\mylabel{Existence}

In this section we establish     existence and uniqueness results for Albanese maps.
These are consequences of a  general criterion, for which we have to generalize   maximal abelian subvarieties
to a relative setting. We start by doing this. Let $S$ be  a base scheme, and 
$G$ be an algebraic space endowed with a group structure. Assume that the structure morphism $G\ra S$ is of finite type.

\begin{definition}
\mylabel{family abelian subvarieties}
A \emph{family   of maximal abelian subvarieties} for $G$ is a family of abelian varieties $A\ra S$,
together with a homomorphism $i: A\ra G$ such that for each point $s\in S$,
the map  identifies the fiber  $A_s$ with the maximal abelian subvariety 
inside the group scheme $G_s$.
\end{definition}

Note that $i:A\ra G$ is a  monomorphism, provided that $G$ is locally separated and $S$ is noetherian,
according to Lemma \ref{kernel trivial}. Then  the translation action of $A$ on $G$ is free,
so  the quotient $G/A$ exists as an algebraic space by Lemma \ref{quotient}.  

Now suppose that $X$ is an algebraic space whose structure morphism $X\ra S$ is  proper, flat, of finite presentation and  cohomologically flat
in degree $d=0$, with $h^0(\O_{X_s})=1$ for all points $s\in S$.  Then  $\Pic^\tau_{X/S}$ is an algebraic space endowed with a group structure
whose structure morphism is  of finite presentation (hence quasiseparated) and locally separated. Note that  it is neither flat nor separated in general.
We can formulate the main result of this paper:

\begin{theorem}
\mylabel{criterion albanese}
Assumptions as above. Then   $G=\Pic^\tau_{X/S}$ admits a   family of maximal abelian subvarieties if and only if 
there is an Albanese map $f:X\ra P$. Moreover, it is universal for morphisms  to families of para-abelian varieties,
and commutes with base-change.
\end{theorem}

Note that in the relative setting, the assumption is restrictive: For example, a Weierstra\ss{} equation
over a discrete valuation ring $R$ whose discriminant $\Delta$ is non-zero and belongs to $\maxid_R$
defines a family $X\subset\PP^2_R$ of cubics where $G=\Pic^\tau_{X/S}$ does not  admit a family of maximal abelian subvarieties:
The generic fiber $G_\eta$ is a one-dimensional abelian variety isomorphic to $X_\eta$, 
whereas the closed fiber of $G$ is isomorphic to the multiplicative or the additive group.
Also note that for the family $X\ra S$ of Enriques surfaces considered in Proposition \ref{enriques counterexample} 
the $\Pic^\tau_{X/S}$ admits a family of maximal abelian subvarieties, namely the zero family.

The proof  of the theorem  is given at the end of this section.
The universal property ensures that the Albanese map is \emph{unique up to unique isomorphism}.
This justifies to write $ \Alb_{X/S}=P$, and we call it  the \emph{family of Albanese varieties}
for $X$. The   Albanese map becomes
$$
f:X\lra \Alb_{X/S}.
$$
This is equivariant  with respect to actions of algebraic spaces endowed with a group structure, an observation  that seems to be new
for infinitesimal actions,  even over ground fields:

\begin{corollary}
\mylabel{albanese equivariant}
Assumptions as in the theorem. Then there is a unique action of  $\Aut_{X/S}$
on the family of Albanese varieties $\Alb_{X/S}$ that makes the Albanese map $f:X\ra\Alb_{X/S}$ equivariant.
\end{corollary}

\proof
Let $\sigma\in\Aut(X)$ be an $S$-automorphism. By the universal property in the theorem, there
is a unique morphism $\sigma_*$ completing the following diagram:
$$
\begin{tikzcd} 
X\ar[r,"\sigma"]\ar[d,"f"']     	     &   X\ar[d,"f"]\\
\Alb_{X/S}\ar[r,dashed,"\sigma_*"']       &   \Alb_{X/S}
\end{tikzcd}
$$
The uniqueness ensures that $\sigma\mapsto\sigma^*$ respects compositions and identities.
In turn,  there is a unique action of the group $\Aut(X)$ making the Albanese map equivariant.

Since the Albanese map commutes with base-change, the same holds for the $\sigma_*$.
Applying the above reasoning with $R$-valued points of $\Aut_{X/S}$ and using the Yoneda Lemma,
we obtain the desired action of the group scheme $\Aut_{X/S}$ making the Albanese map equivariant.
\qed

\medskip
One may reformulate the   theorem in categorical language: Let $\catC$ be the category of all 
algebraic spaces $X$ whose structure morphism $X\ra S$ is proper, flat, of finite presentation,
cohomologically flat in degree $d=0$, with $h^0(\O_{X_s})=1$ for all $s\in S$,
and such that $G=\Pic^\tau_{X/S}$ admits a family of maximal abelian subvarieties.
Write $\catC'\subset \catC$ for the full subcategory comprising all families of para-abelian varieties.
The arguments for the following, which are purely formal and analogous to the preceding proof, are left to the reader:

\begin{corollary}
\mylabel{albanese functorial}
For each morphism $\varphi:X\ra X'$ in the category $\catC$, there is a unique morphism $\varphi_*:\Alb_{X/S}\ra \Alb_{X'/S}$
making the diagram
$$
\begin{CD}
X@>\varphi>> X'\\
@VfVV    @VVf'V\\
\Alb_{X/S}  @>>\varphi_*>    \Alb_{X'/S}
\end{CD}
$$
commutative. Moreover, $\varphi\mapsto\varphi_*$ respects composition and identities, 
and the resulting functor $\catC\ra \catC'$ given by $X\mapsto\Alb_{X/S}$ 
is left adjoint to the inclusion functor $\catC'\ra\catC$.
\end{corollary}

Over ground fields, the assumption in the theorem is vacuous, which gives:

\begin{corollary}
\mylabel{albanese over field}
Suppose that $S $ is the spectrum of  field $k$. Then there is an Albanese map $f:X\ra \Alb_{X/k}$. It is  universal for
morphisms to para-abelian varieties, commutes with field extensions, is equivariant with respect to group scheme actions,
and functorial in $X$.
\end{corollary}

This generalizes previous results for perfect ground fields,  or for geometrically integral schemes,
or geometrically connected and geometrically reduced schemes, compare the discussion in the Introduction.
For results in the non-proper situation,  see appendix in \cite{Wittenberg 2008}, \cite{Achter; Casalaina-Martin; Vial 2022}
and \cite{Schroeer 2022}.
 We can actually formulate an unconditional result in the relative setting, in the spirit of  \cite{SGA 6}, Expos\'e XII, Section 1:

\begin{corollary}
\mylabel{albanese over open}
Suppose that the base scheme $S$ is integral, and that the generic fiber of $\Pic^\tau_{X/S}\ra S$ is proper.
Then after replacing $S$ with some dense open set $U$, there is an  Albanese map $f:X\ra \Alb_{X/S}$. It is  universal for
morphisms into families of para-abelian varieties, commutes with base-change, and is equivariant with respect to actions of relative
group spaces.
\end{corollary}

\proof
Without restriction we may assume that $S$ is the spectrum of an integral noetherian ring $R$.
We have to find some $U$ such that $G=\Pic^\tau_{X/S}$ admits a family of maximal abelian subvarieties over $U$.
Let $\eta\in S$ be the generic point.
Applying Proposition \ref{three step filtration} to the proper group scheme $G_\eta$ over the field of fractions $F=\Frac(R)$,
we see that the maximal abelian variety $A_\eta\subset G_\eta$ is the kernel for the affinization map.

It suffices to treat the case that the algebraic space $G$ 
is schematic, and that $A_\eta$ extends to a family of abelian varieties $A\ra S$,
by shrinking $S$.
In the same way we may assume that the inclusion $A_\eta\subset P_\eta$
extends to some homomorphism $i:A\ra G$, and that the diagonal embedding $G\ra G\times G$ is closed.
Now the kernel  $N=\Kernel(i)$ is a closed subgroup scheme of $A$,
and in particular proper, with $N_\eta=0$. By Chevalley's Semicontinuity Theorem (\cite{EGA IVc}, Corollary 13.1.5)
we may shrink $S$ further making the fibers $N_s$ finite. By generic flatness (\cite{EGA IVb}, Theorem 6.9.1), we even achieve $N=0$.

It remains to verify that  $A_s\subset G_s$ are maximal abelian subvarieties
for a dense open set of points $s\in S$. Let $f:G\ra S$ be the structure morphism, and set $g=\dim(G_\eta)$.
Again by Chevalley's Semicontinuity Theorem (\cite{EGA IVc}, Theorem 13.1.3), the set $Z\subset G$
of all points $x\in G$ with $\dim_x(G_{f(x)})\geq g+1$ is closed. Its image $f(Z)\subset S$ is constructible,
and disjoint from $\eta$. So after shrinking $S$, we may assume that all fibers of $f:G\ra S$ are $g$-dimensional.
For dimension reasons, the inclusion $A_s\subset G_s$ must be the maximal abelian subvariety, for all points $s\in S$.
\qed

\medskip
Moret-Bailly and one referee   alerted us that the conclusion  does not hold without suitable assumption
on the generic fiber, as the following example shows: Let   $X_\QQ$ be the denormalization of an elliptic curve $E_\QQ$ that identifies
the origin $e\in E_\QQ$ with a rational point $a\in E_\QQ$ of  infinite order. Any morphism $X_\QQ\ra P_\QQ$ to an abelian variety
induces a homomorphism $E_\QQ\ra P_\QQ$ with $a$ in the kernel. It follows that these morphisms are constant,
hence $\Alb_{X_\QQ/\QQ}$ is trivial. On the other hand, an extension $X$ over some suitable $R=\ZZ[1/n]$ 
yields a $\Pic^\tau_{X/S}$ that is an extension of a family $E\ra S$ of elliptic curves by the multiplicative group.
Over each closed point $s\in S$, the residue field $\kappa(s)=\FF_p$ is finite, hence the corresponding Ext group is finite, and
the class of the extension has finite order.
It follows that $\Alb_{X_s/\FF_p}$ is one-dimensional.
\bigskip

\emph{Proof of Theorem \ref{criterion albanese}.} 
If there is an Albanese map $f:X\ra P$, the image of the monomorphism $f^*:\Pic^\tau_{P/S}\ra \Pic^\tau_{X/S}$ is a
family maximal abelian subvarieties. Our task is to establish the converse: Suppose there
is a family $A\subset\Pic_{X/S}$ of maximal abelian subvarieties.
We already saw in Proposition \ref{base-change albanese map} that the Albanese map commutes with base-change, once it exists.
Our task here is  to establish   existence and universal property.   We proceed in five intertwined steps,
with various temporary assumptions:
 
\medskip
{\bf Step 1:}
\emph{We show existence, assuming that there is a  Poincar\'e sheaf $\shP$ on the product $X\times\Pic_{X/S}$.} 
Since the family of abelian varieties $A$ has a section, there is also a Poincar\'e sheaf $\shF$ on $ \Pic_{A/S}\times A$.
Tensoring with the preimage of some invertible sheaf on $\Pic^\tau_{A/S}$ we may assume
that $\shF$ becomes trivial on both $\{e\}\times A$ and $\Pic^\tau_{A/S}\times\{e\}$.
Such Poincar\'e sheaves are called \emph{normalized}.
Now regard the restriction $\shP|_{X\times A}$ as a family of invertible sheaves on  
$A$  parameterized by  $X$.
Let $f:X\ra\Pic_{A/S}$ be the classifying map. Then the sheaves
$\shP|_{X\times A}$ and $(f\times\id_A)^*(\shF)$ define the same elements in $\Pic_{A/S}(X)$.
Like \eqref{leray serre} we have an exact sequence
$$
0\lra \Pic(X)\lra\Pic(X\times A)\lra \Pic_{A/S}(X),
$$
so  $\shP|_{X\times A}$ and $(f\times\id_A)^*(\shF)$ differ only by 
the pull-back of some  invertible sheaf on $X$.
Form the relative numerical group  $\Num_{A/S}$ and consider the composite map $\bar{f}:X\ra\Num_{A/S}$.
This   factors over some section $\delta:S\ra\Num_{A/S}$, by Lemma \ref{rigidity}.
The preimage $P=\Pic^\delta_{A/S}$ is a family of para-abelian varieties. It   comes with the structure of a principal homogeneous space
with respect to   $\Pic^\tau_{A/S}$, 
stemming from tensor product of invertible sheaves.
By construction, the classifying map 
factors as  $f:X\ra P$. We claim that this is an Albanese map. 

To see this, it suffices by our very definition
to treat the case that $S$ is the spectrum of an algebraically closed field $k$.
In light of Hilbert's Nullstellensatz, there is a rational point $a\in X$. This induces a rational point
$e=f(a)$ on $P$, and the latter becomes an abelian variety. 
Moreover, the rational point $\delta\in\Num_{A/S}$ comes from an invertible sheaf $\shN$ on $A$.
Tensor products with $\shN$ give  an isomorphism $\Pic^\tau_{A/S}\ra\Pic^\delta_{A/S}$.
Composing $f$ with its inverse, we may assume that $\delta =e$.
The normalized Poincar\'e sheaf $\shF$
induces   identifications $A=\Pic^\tau_{P/k}$ and $P=\Pic^\tau_{A/k}$.
By construction, the abelian varieties $P$ and $A$ have the same dimension, so we merely have to check that
the   kernel $N$ of  $f^*:\Pic^\tau_{P/k}\ra \Pic^\tau_{X/k}$ is trivial.
Suppose this is not the case. Then $N$   contains either a non-zero rational point 
or a tangent vector supported by the origin. In any case, there  is a $k$-algebra $R$ of degree $[R:k]=2$
with an embedding $\Spec(R)\subset N\subset \Pic^\tau_{P/k}=A$
and some non-trivial invertible sheaf $\shL$ on $P\otimes_kR$ such that $(f\otimes\id_R)^*(\shL)$ becomes
trivial on $X\otimes_kR$. In particular, $\shL$ is numerically trivial.
Let $l:\Spec(R)\ra \Pic^\tau_{P/k}=A$ be the classifying map, such that $\shL\simeq \shF|_{P\otimes_kR}$.
It follows that 
$$
(f\otimes\id_R)^*(\shF|_{P\otimes_kR})\simeq (f\otimes\id_R)^*(\shL)\simeq\O_{X\otimes_kR}.
$$
Now recall that $\shP|_{X\times A}$ and $(f\times\id_A)^*(\shF)$ differ by the preimage of some invertible sheaf $\shM$ on $X$.
We infer that $\shP|_{X\otimes_k  R}\simeq \shM\otimes_kR$.
Regarding $\shP|_{X\otimes_kR}$ as a family of invertible sheaves on $X$ parameterized by   $\Spec(R)$,
we conclude that the classifying map $\Spec(R)\ra \Pic_{X/k}$ factors over a closed point.
On the other hand, since $\shP$ is the Poincar\'e sheaf, the classifying map 
is the composition of the   embeddings $\Spec(R)\subset A$ and $A\subset\Pic_{X/k}$,
contradiction.

\medskip
{\bf Step 2:}
\emph{We verify the universal property for the particular Albanese map above,
assuming that $X$ and all the algebraic spaces  $\Pic^\delta_{X/S}$, $\delta\in\Num_{X/S}(S)$ admit sections.}
Since $X$ has a section, there is   a  Poincar\'e sheaf $\shP$ on $X\times\Pic_{X/S}$. Let $f:X\ra P$ be the resulting Albanese map constructed in step 1,
and $g:X\ra Q$ be a morphism into another family of para-abelian varieties. We have to verify that $g$ factors over $f$, via some  $h:P\ra Q$.
Note that we already saw in Proposition \ref{weak uniqueness} that such a factorization is unique, once it exists.
Fix a section $s:S\ra X$. This induces a section for $Q\ra S$, which therefore becomes
a family of abelian varieties. Choose a normalized Poincar\'e sheaf $\shG$ on $\Pic^\tau_{Q/S}\times Q$.
Viewing this as a family of invertible sheaves on $\Pic^\tau_{Q/S}$ parameterized by $Q$, we see that
 $g:X\ra Q$ is the classifying map for each of the sheaves
\begin{equation}
\label{pullback poincare}
\shM=(\id\times g)^*(\shG)\otimes\pr_X^*(\shL)
\end{equation}
on $\Pic^\tau_{Q/S}\times X$, where $\shL$ is any invertible sheaf on $X$.
Our task is to find  $\shL$ such that   $\shM\simeq (\id\times f)^*(\shN)$ for some
invertible sheaf $\shN$ on $\Pic^\tau_{Q/S}\times P$, which then gives  the desired factorization $g=h\circ f$. 
We shall achieve this  by successively constructing the
dashed arrows in the following commutative diagram, starting with $\varphi\times\id$ and proceeding clockwise.
For increased clarity, we write $\Pic_Q$ rather than $\Pic_{Q/S}$  etc.:
\begin{equation}
\label{diamond diagram}
\begin{tikzcd}[]
                            & \Pic^\tau_Q\times P\arrow[dl,dashed,"\psi\times \id"']\arrow[rr,dashed,"\id\times h"]      &                           & \Pic^\tau_Q\times Q\\
\Pic^\tau_P\times P     &                               & \Pic^\tau_Q\times X\arrow[ur,"\id\times g"]\arrow[dr,"g^*\times \id"]\arrow[ul,"\id\times f"']\arrow[dd,dashed,"\varphi\times\id"]\arrow[dl,dashed,"\psi\times\id"']\\
                            & \Pic^\tau_P\times X\arrow[ul,"\id\times f"]\ar[dr,"f^*\times\id"']         &  & \Pic_X\times X \\
                            &                                                                               & \Pic^\tau_X\times X \ar[ru,hook,"\can\times\id"']
\end{tikzcd}	 
\end{equation} 
Fix any invertible sheaf $\shL$ on $X$, 
regard the resulting sheaf $\shM$ in \eqref{pullback poincare}   as a family of invertible sheaves on $X$ parameterized by $\Pic^\tau_{Q/S}$, and consider the   
classifying map $\varphi:\Pic^\tau_{Q/S}\ra\Pic_{X/S}$.
In light of   Lemma \ref{rigidity}, the composite map $\overline{\varphi}:\Pic^\tau_{Q/S}\ra\Num_{X/S}$ factors over some section $\delta:S\ra\Num_{X/S}$.
By assumption  $P=\Pic^\delta_{X/S}$ admits a section. The latter comes from an invertible sheaf $\shL_0$ on $X$,
because   $X$   has a section, too. Replacing $\shL$ with $\shL\otimes\shL_0^{\otimes-1}$, we may assume that $\delta$ is the zero section.
This yields  a factorization
$\varphi:\Pic^\tau_{Q/S}\ra \Pic^\tau_{X/S}$ and gives the vertical dashed arrow in \eqref{diamond diagram}. 
By Proposition \ref{weak uniqueness}, it comes with a factorization
 $\psi:\Pic^\tau_{Q/S}\ra\Pic^\tau_{P/S}$, yielding the two  diagonal dashed arrows.
By construction, \eqref{pullback poincare} is the pullback of the normalized Poincar\'e sheaf $\shF$ on $\Pic^\tau_{P/S}\times P$.
Hence it is the pullback of $\shN=(\psi\times\id_P)^*(\shF)$, giving the desired factorization $g=h\circ f$ of classifying maps
and the upper dashed arrow in \eqref{diamond diagram}.

\medskip
{\bf Step 3:}
\emph{We verify the universal property for the particular Albanese map above,
now only assuming that $X$  admits a section.} Keep the notation from the previous step.
There the arguments only relied on the existence of   sections in $\Pic^\delta_{X/S}$ for one particular
$\delta\in \Num_{X/S}(S)$ occurring along the way. Choose an fppf morphism $S'\ra S$ so that
this $\Pic^\delta_{X/S}$ acquires an $S'$-valued point.
Then step 2 gives the desired factorization $g'=h'\circ f'$ for the base-changes $X'=X\times_SS'$, $P'=P\times_SS'$
and $Q'=Q\times_SS'$. The two pullbacks of $h'$ to $S'\times_SS'$ coincide, by the uniqueness in Proposition \ref{weak uniqueness}.
Now fppf descent (\cite{SGA 1}, Expos\'e VIII, Theorem 5.2 for schemes and \cite{SP}, 
Lemma 0ADV for algebraic spaces) 
gives the desired factorization $h:P\ra Q$.

\medskip
{\bf Step 4:}
\emph{We prove the universal property for general Albanese maps $f:X\ra P$.}
Let $g:X\ra Q$ be another morphism into some family of para-abelian varieties. As above we have  to find a factorization
$g=h\circ f$, and we shall achieve this by descent. Note that   uniqueness was already established in Proposition \ref{weak uniqueness}.
Choose some fppf morphism $S'\ra S$ so that the base-change $X'=X\times_SS'$ admits a section.
Set $P'=P\times_SS'$, and let $X'\ra P_1'$ be the Albanese map over $S'$ constructed in step 1 with Poincar\'e sheaves.
We then have a factorization $P_1'\ra P'$, according to step 3. By our definition of Albanese maps,
the induced proper homomorphism $\Pic^\tau_{P_1'/S'}\ra\Pic^\tau_{P'/S'}$ is fiberwise   an isomorphism, 
thus a monomorphism (Lemma \ref{kernel trivial}), hence a bijective closed embedding (\cite{EGA IVd}, Corollary 18.12.6).
With the Nakayama Lemma one infers that $\Pic^\tau_{P_1'/S'}\ra\Pic^\tau_{P'/S'}$ is an isomorphism.
By biduality, the original homomorphism $P'_1\ra P'$ is an isomorphism.
Regard it as an identification $P'=P'_1$. Using step 3 again, we get a factorization 
$g'=h'\circ f'$.
The two pullbacks of $h':P'\ra Q'$ to $S''=S'\times_SS'$ coincide, by   uniqueness.
Now \cite{SP},  
Lemma 0ADV gives the desired factorization $h:P\ra Q$.

\medskip
{\bf Step 5:}
\emph{We establish existence in general.}
Choose some fppf morphism $S'\ra S$ so that the base-change $X'=X\times_SS'$ acquires a section.
This happens, for example, with $S'=X$.
Now step 1 gives an Albanese map $f':X'\ra P'$.  Set $S''=S'\times_SS'$, and consider the two base-changes
$\pr_1^*(P')$ and $\pr_2^*(P')$. The  two induced morphisms from $X'=X\times_SS''$ are related by a  morphism $\varphi:\pr_1^*(P')\ra\pr_2^*(P')$,
according to the existence part in step 4, applied over $S''$, and uniqueness ensures that $\varphi$ is an isomorphism.
This satisfies the cocycle condition over  $S'\times_SS'\times_SS'$, in light of the uniqueness part in step 4.
In turn, the scheme $P'$ over $S'$ descends to an algebraic space $P$ over $S'$.
By construction, $P'\ra S'$ is a family of abelian varieties, hence $P\ra S$ is a family of para-abelian varieties.
In the same way, the Albanese map $f':X'\ra P'$ descends to a morphism $f:X\ra P$. The latter
is an Albanese map, by Proposition \ref{base-change albanese map}.
\qed

\section{Appendix:  Embeddings for algebraic spaces}
\mylabel{Appendix}
 
Here we collect some general facts about embeddings of algebraic spaces that were used throughout, and seem to be of independent
interest. The topic indeed requires   attention,
because   schemes are locally separated, whereas not all algebraic spaces share this property.
Recall that a morphism of schemes $f:Z\ra Y$ is   an \emph{embedding} 
if it factors into a closed embedding $Z\ra U$ followed by an open embedding $U\ra Y$.
Then the set-theoretical image $C=f(Z)$ is locally closed, and $Z$ is determined, up 
to unique isomorphism,   by  an open set $U\subset Y$ and a quasicoherent     $\shI\subset \O_U$.
Note that such   $U$ are not unique, but there is a maximal one, namely
the complement of  the closed set $\bar{C}\smallsetminus C=\bar{C}\cap (Y\smallsetminus U)$.

A morphism of algebraic spaces $f:Z\ra Y$ is called an \emph{embedding} if for
each affine \'etale neighborhood $V\ra Y$ the fiber product $Z\times_YV$ is a scheme,
and the projection $Z\times_YV\ra V$ is an embedding.  

\begin{lemma}
\mylabel{embedding spaces}
Let  $f:Z\ra Y$ be an embedding of algebraic spaces. Suppose that the morphism is quasicompact.
Then it admits  a factorization into a closed embedding $Z\ra U$ followed
by an open embedding $U\ra Y$.
\end{lemma}

\proof
Let $V_\lambda \ra Y$, $\lambda\in I$ be the affine \'etale neighborhoods, and $C_\lambda \subset V_\lambda $ 
be the set-theoretical image  of the projection $Z\times_YV_\lambda\ra V_\lambda$. This is locally closed, 
so  $\bar{C}_\lambda\smallsetminus C_\lambda$ is a closed set in $V_\lambda$. Moreover,  the scheme $Z\times_YV_\lambda$
is quasicompact, by the corresponding assumption on $f$.
For each  morphism of neighborhoods $i:V_\lambda\ra V_\mu$    we have $i^{-1}(C_\mu)=C_\lambda$.
Since the morphism  is \'etale and hence flat, and the scheme $Z\times_YV_\lambda$ is quasicompact,
we also have $i^{-1}(\bar{C}_\mu) = \bar{C}_\lambda$ by \cite{SGA 1}, Expos\'e VIII, Theorem 4.1, and thus
$i^{-1}(\bar{C}_\mu\smallsetminus C_\mu)=\bar{C}_\lambda\smallsetminus C_\lambda$.
Thus the complementary  open sets $U_\lambda\subset V_\lambda$ are compatible, and thus define an open embedding $U\ra Y$,
according to \cite{SP}, Lemma 0ADV. Note that in the latter result, we indeed can ignore the condition on the cardinality of the index set $I$   by working 
in a fixed Grothendieck universe, confer the discussion in \cite{Schroeer 2017}, Section 2.
By construction, the morphism $Z\ra Y$ factors over $U$, and the resulting $Z\ra U$ is a closed embedding.
\qed

\begin{proposition}
\mylabel{embedding yoga}
Let $f:X\ra Y$ and $g:Y\ra Z$ be morphisms of algebraic spaces. If the composition $g\circ f$ is an embedding,
and $g$ is locally separated, then $f$ is an embedding.
\end{proposition}

\proof
This follows as in \cite{EGA I}, Section 5.2. We recall the argument for convenience:
In the diagram
$$
\begin{tikzcd} 
    & Y\times_ZY                                    &               & Z\times_ZY    & \\
Y\ar[ur,"\Delta"]  &                                               & X\times_ZY\ar[ul,"f\times\id"']\ar[ur,"gf\times\id"]    &               &  Y\ar[ul,"\simeq"']   \\
    & X\times_YY\ar[rr,"\simeq"']\ar[ul]\ar[ur,"\id\times\id"]\ar[ul,"f\times\id"]       &               & X\ar[ur,"f"']\ar[ul,"\Gamma_f"']\\   
\end{tikzcd}
$$
the square to the left is cartesian. Since the arrows $\Delta$ and $gf$ are embeddings, the same holds for 
$\id\times \id$ and $gf\times\id$, and thus also for $f$.
\qed

\medskip
The following application  was a crucial technical step in the proof for Lemma \ref{rigidity}:

\begin{corollary}
\mylabel{section embedding}
Let   $N$ be a locally separated algebraic space. For each $\sigma\in N(S)$
the morphism $\sigma:S\ra N$ is an embedding.
\end{corollary}

\proof
The structure morphism  $\varphi:N\ra S$ is locally separated, and the composition
$\varphi\circ\sigma=\id_S$ is an embedding. By the Proposition, also $\sigma$ is an embedding.
\qed

\medskip
Let us record the following useful consequence for homomorphisms $f:G\ra N$ between algebraic spaces $G,N$ 
endowed with    group structures:

\begin{lemma}
\mylabel{kernel trivial}
In the above situation, suppose that  the structure morphism $N\ra S$ is locally separated, and that $G\ra S$ is separated and of finite type.
If  $ G_s\ra N_s$ has trivial kernel for all points $s\in S$, then $f:G\ra N$ is a monomorphism.
\end{lemma}

\proof
We have to show that  $H=\Kernel(f)$ is trivial. This kernel is given by a fiber product $H=G\times_N\{e_S\}$.
By Corollary \ref{section embedding}, the neutral section $e:S\ra N$ is an embedding, so the same holds for the base-change $H\ra G$.
It follows that the structure morphism $H\ra S$ is separated and of finite type.
Obviously, the formation of kernels commutes with base-change, and the fibers $H_s$ are trivial. Thus $H\ra S$ is universally bijective,
and in particular quasi-finite. This ensures that the algebraic space  $H$ is schematic (\cite{Laumon; Moret-Bailly 2000}, Theorem A.2).
The relative group scheme $H\ra S$ is thus trivial, by \cite{SGA 3b}, Expos\'e VI$_B$,  Corollary 2.10.
\qed


\end{document}